\documentclass[aop,MSNbibl,seceqn,dvips]{arximspdf}

\doi{10.1214/11-AOP665}
\volume{40}
\issue{4}
\pubyear{2012}
\firstpage{1829}
\lastpage{1859}

\makeatletter

\newcommand{\iint}{\int\!\!\!\int}
\newcommand{\iiint}{\int\!\!\!\int\!\!\!\int}

\let\epsilon\varepsilon
\newcommand{\dimh}{\dim_{\mathrm{H}}}
\newcommand{\eps}{\epsilon}
\newcommand{\1}{\mathbf{1}}
\renewcommand{\d}{\mathrm{d}}
\renewcommand{\P}{\mathrm{P}}
\newcommand{\E}{\mathrm{E}}
\newcommand{\R}{\mathbf{R}}
\newcommand{\F}{\mathcal{F}}
\newcommand{\eqref}[1]{(\ref{#1})}

\newtheorem{theorem}{Theorem}[section]
\newtheorem{corollary}[theorem]{Corollary}
\newtheorem{proposition}[theorem]{Proposition}
\newtheorem{lemma}[theorem]{Lemma}

\newproclaim{definition}[theorem]{Definition}
\newproclaim{remark}[theorem]{Remark}

\makeatother

\begin{document}
\begin{frontmatter}

\title{Critical Brownian sheet does not have double~points}
\runtitle{Double points of Brownian sheet}

\begin{aug}
\author[A]{\fnms{Robert C.} \snm{Dalang}\ead[label=e1]{robert.dalang@epfl.ch}\ead[url,label=u1]{http://mathaa.epfl.ch/\textasciitilde rdalang}\thanksref{aut1}},
\author[B]{\fnms{Davar} \snm{Khoshnevisan}\corref{}\ead[label=e2]{davar@math.utah.edu}\ead[url,label=u2]{http://www.math.utah.edu/\textasciitilde davar}\thanksref{aut2}},
\author[C]{\fnms{Eulalia} \snm{Nualart}\ead[label=e3]{nualart@math.univ-paris13.fr}\ead[url,label=u3]{http://www.nualart.es}},
\author[D]{\fnms{Dongsheng} \snm{Wu}\ead[label=e4]{Dongsheng.Wu@uah.edu}\ead[url,label=u4]{http://webpages.uah.edu/\textasciitilde dw0001}}
and
\author[E]{\fnms{Yimin} \snm{Xiao}\ead[label=e5]{xiao@stt.msu.edu}\ead[url,label=u5]{http://www.stt.msu.edu/\textasciitilde xiaoyimi}\thanksref{aut2}}
\runauthor{R. C. Dalang et al.}
\thankstext{aut1}{Supported by grant from the Swiss National Science Foundation.}
\thankstext{aut2}{Supported by grant from the US National Science Foundation.}
\affiliation{EPF-Lausanne, University of Utah, University of Paris
13, University of Alabama--Huntsville and Michigan State University}
\address[A]{R. C. Dalang\\
Institut de Math\'ematiques\\
Ecole Polytechnique
F\'ed\'erale de Lausanne\\
Station 8\\
CH-1015 Lausanne\\
Switzerland\\
\printead{e1}\\
\printead{u1}} 
\address[B]{D. Khoshnevisan\\
Department of Mathematics\\
University of Utah\\
Salt Lake City, Utah 84112--0090\\
USA\\
\printead{e2}\\
\printead{u2}\hspace*{24pt}}
\address[C]{E. Nualart\\
Institut Galil\'ee\\
Universit\'e Paris 13\\
93430 Villetaneuse\\
France\\
\printead{e3}\\
\printead{u3}}
\address[D]{D. Wu\\
Department of Mathematical Sciences\\
University of Alabama--Huntsville\\
Huntsville, Alabama 35899\\
USA\\
\printead{e4}\\
\printead{u4}}
\address[E]{Y. Xiao\\
Department of Statistics and Probability\\
Michigan State University\\
East Lansing, Michigan 48824\\
USA\\
\printead{e5}\\
\printead{u5}}
\end{aug}

\received{\smonth{9} \syear{2010}}
\revised{\smonth{3} \syear{2011}}

%
\begin{abstract}
We derive a decoupling formula for the Brownian sheet
which has the following ready consequence: An $N$-parameter
Brownian sheet in~$\mathbf{R}^d$ has double points if and only
if $d<4N$. In particular, in the critical case where
$d=4N$, the Brownian sheet does not have double points.
This answers an old problem in the folklore of the subject.
We also discuss some of the geometric consequences
of the mentioned decoupling, and establish a partial
result concerning $k$-multiple points in the critical
case $k(d-2N) = d$.
\end{abstract}

%

\end{frontmatter}

\section{Introduction}\label{sec1}

Let $B:=(B^1,\ldots,B^d)$ denote a $d$-dimensional
$N$-para\-meter Brownian sheet.
That is, $B$ is a $d$-dimensional,
$N$-parameter, centered Gaussian process
with
%
\begin{equation}
\operatorname{Cov}(B^i(\mathbf{s}) ,B^j(\mathbf{t}))=\delta
_{i,j}\cdot
\prod_{k=1}^N(s_k\wedge t_k),
\end{equation}
where $\delta_{i,j}=1$ if $i=j$ and $0$ otherwise,
and $\mathbf{s},\mathbf{t}\in\R^N_+$, $\mathbf{s} = (s_1,\dots,s_N)$,
$\mathbf{t} = (t_1,\dots,t_N)$. Here and throughout, we define
%
\begin{equation}
\bolds{\mathcal{T}}:= \{ (\mathbf{s} ,\mathbf{t})\in(0,\infty
)^{2N}\dvtx
s_i \neq t_i\mbox{ for all }i=1,\dots,N \}.
\end{equation}
The following is the main result of this paper.

\begin{theorem}\label{thmain}
Choose and fix a Borel set $A\subseteq\R^d$. Then
%
\begin{equation}
\P\{ \exists (\mathbf{u}_1,\mathbf{u}_2)\in\bolds{\mathcal{T}}\dvtx
B(\mathbf{u}_1)=B(\mathbf{u}_2)\in A \}>0
\end{equation}
if and only if
%
\begin{equation}
\P\{ \exists (\mathbf{u}_1,\mathbf{u}_2)\in\bolds{\mathcal{T}}\dvtx
W_1(\mathbf{u}_1)=W_2(\mathbf{u}_2)\in A \}>0,
\end{equation}
where $W_1$ and $W_2$ are independent $N$-parameter
Brownian sheets in $\R^d$ (unrelated to
$B$).
\end{theorem}

Theorem \ref{thmain} helps answer various questions
about the multiplicities of the random surface
generated by the Brownian sheet.
We introduce some notation in order to present
some of these issues.

Recall that $x\in\R^d$ is a
\emph{$k$-multiple point of $B$}
if there exist distinct points $\mathbf{s}_1,\ldots,\mathbf{s}_k\in
(0,\infty)^N$ such that $B(\mathbf{s}_1)=\cdots=B(\mathbf{s}_k)=x$.
We write $M_k$ for the collection of
all $k$-multiple points of $B$.
Note that $M_{k+1} \subseteq M_k$ for all $k \ge2$.

In this paper, we are concerned mainly with the case $k=2$;
elements of~$M_2$ are the \emph{double points} of $B$.
In Section \ref{secpfmain} below, we
derive the following ready consequence of Theorem \ref{thmain}.

\begin{corollary}\label{cohit}
Let $A$ denote a nonrandom Borel set in $\R^d$.
If $d>2N$, then
%
\begin{equation}
\P\{ M_2 \cap A\neq\varnothing\}>0\quad
\mbox{if and only if}\quad
\operatorname{Cap}_{2(d-2N)}(A)>0,
\end{equation}
where $\operatorname{Cap}_\beta$ denotes the Bessel--Riesz
capacity in dimension $\beta\in\R$; see Section \ref{seccap}
below. If $d=2N$, then
$\P\{ M_2 \cap A\neq\varnothing\}>0$
if and only if there exists a probability measure
$\mu$, compactly supported in $A$, such that
%
\begin{equation}\label{eqbilog}
\iint \biggl| \log_+ \biggl(\frac{1}{|x-y|} \biggr)
\biggr|^2 \mu(\d x) \mu(\d y)<\infty.
\end{equation}
Finally, if $d<2N$, then $\P\{M_k\cap A\neq\varnothing\}>0$
for all $k \ge2$ and all nonvoid, nonrandom Borel sets $A\subset\R^d$.
\end{corollary}

We apply Corollary \ref{cohit} with $A:=\R^d$ and appeal
to Taylor's theorem (\cite{Khbook}, pages 523--525), to deduce
the following.

\begin{corollary}\label{coMP}
An $N$-parameter, $d$-dimensional
Brownian sheet has double points if and only if $d<4N$.
In addition, $M_2$ has positive Lebesgue
measure almost surely if and only if $d<2N$.
\end{corollary}

When $N=1$, $B$ is $d$-dimensional
Brownian motion, and this corollary has a rich history in that case:
L\'evy \cite{Levy} was the first to prove that Brownian motion has double
points ($M_2\neq\varnothing$) when $d= 2$; this is also
true in one dimension, but almost tautologically so. Subsequently,
Kakutani \cite{Kakutani}
proved that Brownian motion in $\R^d$ does not have double points when
$d\ge5$;
see also Ville \cite{Ville}.
Dvoretzky, Erd\H{o}s and Kakutani \cite{DEK50}
then showed that Brownian motion has double points when $d=3$,
but does not have double points in the case\vadjust{\goodbreak} that $d=4$.
Later on, Dvoretzky, Erd\H{o}s and Kakutani \cite{DEK54} proved that
in fact, $M_k\neq\varnothing$ for all $k\ge2$,
when $d=2$.
The remaining case is that $M_3\neq\varnothing$ if and only if
$d\le2$; this fact is due to Dvoretzky et al.
\cite{DEKT}.

When $N>1$ and $k=2$, Corollary \ref{coMP} is new
only in the critical case where
$d=4N$. The remaining (noncritical)
cases are much simpler to derive
and were worked out earlier by one of us \cite{Kpolar}.
In the critical case, Corollary~\ref{coMP} asserts that Brownian sheet has no double
points. This justifies the title of the paper and solves
an old problem in the folklore of the subject.
For an explicit mention---in print---of this problem,
see Fristedt's
review of the article of Chen
\cite{Chen} in \emph{Mathematical Reviews}, where
most of the assertion about $M_2$ (and even
$M_k$) having positive measure was conjectured.

The proof of Theorem \ref{thmain} leads to another interesting
property, whose description requires us first to introduce
some notation.
We identify subsets of $\{1,\ldots,N\}$ with
partial orders on $\R^N$ as follows \cite{KX}:
For all $\mathbf{s},\mathbf{t}\in\R^N$ and $\pi\subseteq\{
1,\ldots,N\}$,
%
\begin{equation}
\mathbf{s}\prec_{\pi}\mathbf{t}
\quad \mbox{iff}\quad
\cases{
s_i\le t_i,&\quad for all $i\in\pi$,\cr
s_i\ge t_i,&\quad for all $i\notin\pi$.
}
\end{equation}
Clearly every $\mathbf{s}$ and $\mathbf{t}$ in $\R^N$ can be
compared via some $\pi$. In fact, $\mathbf{s}\prec_\pi\mathbf{t}$,
where $\pi$ is the collection of all $i\in\{1,\ldots,N\}$
such that $s_i\le t_i$. We might write $\mathbf{s}\prec_\pi\mathbf{t}$
and $\mathbf{t}\succ_\pi\mathbf{s}$ interchangeably. Sometimes, we will
also write $\mathbf{s}\curlywedge_\pi\mathbf{t}$
for the $N$-vector whose $j$th coordinate is $\min(s_j,t_j)$
if $j\in\pi$ and $\max(s_j,t_j)$ otherwise.

Given a partial order $\pi\subset\{1,\dots,N\}$ and
$\mathbf{s}, \mathbf{t} \in(0,\infty)^N$ we write
$\mathbf{s} \ll_\pi\mathbf{t}$ if $\mathbf{s} \prec_{\pi}
\mathbf{t}$ and
$s_i \neq t_i$, for all $i \in\{1,\dots,N\}$. Define
%
\begin{equation}
\widetilde{M}_k := \left\{
x \in\R^d \left|
\matrix{\exists \mathbf{s}_1,\dots,\mathbf{s}_k \in
(0,\infty)^N\dvtx
B(\mathbf{s}_1)=\cdots=B(\mathbf{s}_k)=x\cr
\mbox{and }
\mathbf{s}_1 \ll_\pi\cdots\ll_\pi\mathbf{s}_k \mbox{ for some }
\pi\subset\{1,\dots,N\}
}
\right.\right\}.\hspace*{-35pt}
\end{equation}

\begin{proposition}\label{comonotone}
Let $A\subset\R^d$ be a nonrandom Borel set. Then for all $k\ge2$,
%
\begin{equation}
\P\{\widetilde{M}_k \cap A \neq\varnothing\} >0
\quad \mbox{if and only if}\quad  \operatorname{Cap}_{k(d-2N)}(A)>0.
\end{equation}
In particular, there are (strictly) $\pi$-ordered $k$-tuples on which
$B$ takes a~common value if and only if $k(d-2N)<d$.
\end{proposition}

Theorem \ref{thmain} can also be used to study various
geometric properties of the random set $M_2$ of double points of $B$.
Of course, we need to study only the case where
$M_2\neq\varnothing$ almost surely. That is, we assume henceforth
that $d<4N$. With this convention in mind,
let us start with the following formula:\looseness=-1
%
\begin{equation}\label{eqdimhM}
\dimh M_2 = d-2(d-2N)^+ \qquad \mbox{almost surely.}
\end{equation}\looseness=0
This formula appears in Chen \cite{Chen} (with a gap in his proof
that was filled by Khoshnevisan, Wu and Xiao \cite{KWX}). In fact, a formula
for $\dimh M_k$ analogous to (\ref{eqdimhM}) holds for all $k \ge2$
\cite{Chen,KWX} and has many connections to the well-known
results of Orey and Pruitt \cite{OP}, Mountford \cite{Mountford}
and Rosen \cite{Rosen}.\vadjust{\goodbreak}

As yet another application of Theorem \ref{thmain} we can refine
\eqref{eqdimhM}
by determining the Hausdorff dimension of $M_2\cap A$ for any nonrandom
closed set
$A\subset\R^d$. First, let us remark that a standard covering
argument [similar to
the proof of part (i) of Lemma \ref{LemExtra}] shows that for any
fixed nonrandom
Borel set $A \subset\R^d$
%
\begin{equation}\label{EqM2Adim}
\dimh(M_2 \cap A) \le\dimh A - 2 (d-2N) \qquad \mbox{almost surely.}
\end{equation}
The following corollary provides an essential lower bound for \mbox{$\dimh (M_2 \cap A)$}.
Recall that the essential supremum $\|Z\|_{L^\infty(\P)}$ of
a nonnegative random variable $Z$ is defined as
%
\begin{equation}
\|Z\|_{L^\infty(\P)} := \inf\bigl\{\lambda>0\dvtx
\P\{Z>\lambda\}=0 \bigr\}\qquad
(\inf\varnothing:=+\infty).
\end{equation}

\begin{corollary}\label{comain}
Choose and fix a nonrandom closed set $A\subset\R^d$. If $\dimh A< 2(d-2N)$,
then with probability one $A$ does not contain
any double points of the Brownian sheet.
On the other hand, if $\dimh A\ge2(d-2N)$, then
%
\begin{equation}\label{eqdimh}
\| \dimh(M_2 \cap A) \|_{L^\infty(\P)}
= \dimh A - 2(d-2N)^+.
\end{equation}
\end{corollary}

Equation \eqref{eqdimhM} follows from Corollary \ref{comain}
and the fact that $\dimh M_2$ is a.s. a~constant. The proof of this
``zero--one law'' follows more-or-less standard methods, which we skip.

There is a rich literature of decoupling, wherein expectation
functionals for sums of dependent random variables are
analyzed by making clever comparisons to similar
expectation functionals that involve only sums of
\emph{independent} (sometimes conditionally independent)
random variables. For a~definitive account,
see the recent book of de la Pe\~{n}a and Gin\'e
\cite{delaPenaGine}.

Theorem \ref{thmain} of the present paper
follows the general philosophy of
decoupling, but applies it to random fractals rather than
random variables (or vectors). A ``one-parameter'' version
of these ideas appears earlier in the work of Peres \cite{Peres96}.
From a technical point of view, Theorem \ref{thmain} is rather
different from the results of decoupling theory.

This paper is organized as follows. Section \ref{seccap}
recalls the main notions of potential theory and presents our
main technical result concerning conditional laws of the Brownian
sheet (Theorem \ref{thCPT}). In Section \ref{sec3}, we present a
sequence of estimates concerning the pinned Brownian sheet. Section
\ref{sec4} contains the proof of Theorem \ref{thCPT}. Finally,
Section \ref{secpfmain} contains the proofs of Theorem \ref{thmain},
of its corollaries and of Proposition \ref{comonotone}.

\section{Potential theory}
\label{seccap}

In this section, we first introduce
some notation for capacities, energies and Hausdorff
dimension, and we also recall some
basic facts about them.
Then we introduce the main technical result of
this paper, which is a theorem of ``conditional potential theory''
and is of independent interest.

\subsection{Capacity, energy, and dimension}
For all real numbers $\beta$, we define a function
$\kappa_\beta\dvtx\R^d\to\R_+\cup\{\infty\}$ as follows:
%
\begin{equation}\label{eqkappa}
\kappa_\beta(x) :=
\cases{
\|x\|^{-\beta},&\quad  if $\beta>0$,\cr
\log_+(\|x\|^{-1}),&\quad if $\beta=0$,\cr
1,&\quad  if $\beta<0$,
}
\end{equation}
where, as usual, $1/0:=\infty$ and $\log_+(z):=1\vee\log(z)$
for all $z\ge0$.

Let $\mathcal{P}(G)$ denote the collection of all probability
measures that are supported by the Borel set $G\subseteq\R^d$,
and define the $\beta$-dimensional \emph{capacity} of $G$ as
%
\begin{equation}
\operatorname{Cap}_\beta(G) := \Bigl[\mathop{\inf_{\mu\in\mathcal{P}(K):}}_{K\subset G\mathrm{\ is\ compact}}
\mathrm{I}_\beta(\mu) \Bigr]^{-1},
\end{equation}
where $\inf\varnothing:=\infty$, and $\mathrm{I}_\beta(\mu)$ is
the $\beta$-dimensional \emph{energy} of $\mu$, defined
as follows, for all $\mu\in\mathcal{P}(\R^d)$ and $\beta\in\R$:
%
\begin{equation}
\mathrm{I}_\beta(\mu) := \iint\kappa_\beta(x-y) \mu(\d x)
\mu(\d y).
\end{equation}
In the cases where $\mu(\d x)=f(x)\,\d x$, we may also write
$\mathrm{I}_\beta(f)$ in place of~$\mathrm{I}_\beta(\mu)$.

Let us emphasize that for all probability measures $\mu$ on $\R^d$
and all Borel sets $G\subseteq\R^d$,
%
\begin{equation}
\mathrm{I}_\beta(\mu)=
\operatorname{Cap}_\beta(G)= 1
\qquad \mbox{when }\beta<0.
\end{equation}

According to Frostman's theorem (\cite{Khbook}, page 521),
the Hausdorff dimension of $G$ satisfies
\begin{eqnarray}\label{eqfrostman}
\dimh G &=& \sup\{\beta>0\dvtx \operatorname{Cap}_\beta(G)>0 \}\nonumber\\[-8pt]\\[-8pt]
&=&\inf\{\beta>0\dvtx \operatorname{Cap}_\beta(G)=0 \}.\nonumber
\end{eqnarray}
The reader who is unfamiliar with Hausdorff dimension
can use the preceding as its \emph{definition}.
The usual definition can be found
in Appendix C of Khoshnevisan \cite{Khbook}, where many properties
of $\dimh$ are also derived. We will also need the following
property:
%
\begin{equation}\label{eqtaylor}
\operatorname{Cap}_n(\R^n)=0 \qquad \mbox{for all }n\ge1.
\end{equation}
See Corollary 2.3.1 of Khoshnevisan (\cite{Khbook}, page 525) for
a proof.

\subsection{Conditional potential theory}

Throughout, we assume that our underlying probability
space $(\Omega,\F,\P)$ is complete.
Given a partial order $\pi$ and a point $\mathbf{s}\in\R^N_+$,
we define $\F_\pi(\mathbf{s})$ to be the $\sigma$-algebra
generated by $\{B(\mathbf{u}),  \mathbf{u}\prec_\pi\mathbf{s}\}$
and all $\P$-null sets. We 
then make the filtration $(\F_\pi(\mathbf{s}), s \in\R^N_+)$
right-continuous in the partial order $\pi$, so that
$\F_\pi(\mathbf{s})=\bigcap_{\mathbf{t}\succ_\pi\mathbf{s}} \F
_\pi(\mathbf{t})$.

\begin{definition}
Given a sub-$\sigma$-algebra $\mathcal{G}$ of
$\F$ and a set-valued function~$A$---mapping $\Omega$ into\vadjust{\goodbreak} subsets of $\R^d$---we say that $A$
is a \emph{$\mathcal{G}$-measurable random set} if
$\Omega\times\R^d\ni
(\omega,x)\mapsto\1_{A(\omega)}(x)$ is $(\mathcal{G}\times
\mathcal{B}(\R^d))$-measurable, where~$\mathcal{B}(\R^d)$
denotes the Borel $\sigma$-algebra on $\R^d$.
\end{definition}

We are also interested in two variants of this definition. The
first follows:

\begin{definition}
Given a $\sigma$-algebra $\mathcal{G}$ of
$\F$, we
say that $f\dvtx\Omega\times\R^d\to\R_+$ is a \emph{%
$\mathcal{G}$-measurable random probability density function} when
$f$ is $(\mathcal{G}\times
\mathcal{B}(\R^d))$-measurable and
$\P\{ \int_{\R^d}f(x)\,\d x=1\}=1$.
\end{definition}

The second variant is the following:
\begin{definition}
Given a $\sigma$-algebra $\mathcal{G}$ of
$\F$, we say that $\rho\dvtx\Omega\times\mathcal{B}(\R^d)\to
[0,1]$ is a \emph{%
$\mathcal{G}$-measurable random probability measure} when
both of the following hold:
\begin{enumerate}
\item$\Omega\ni\omega\mapsto\rho(\omega,A)$ is
$\mathcal{G}$-measurable for every $A \in\mathcal{B}(\R^d)$;
\item$A\mapsto\rho(\omega,A)$ is a Borel probability measure
on $\R^d$ for almost every $\omega\in\Omega$.
\end{enumerate}
\end{definition}

For all $\pi\subseteq\{1,\ldots,N\}$ and
$\mathbf{s}\in\R^N_+$, let $\P^\pi_{\mathbf{s}}$ be a
regular conditional distribution for $B$ given $\F_\pi(\mathbf{s})$,
with the corresponding expectation operator written as $\E^\pi
_{\mathbf{s}}$. That is,
%
\begin{equation}
\E^\pi_{\mathbf{s}} f := \int f\, \d\P^\pi_{\mathbf{s}} = \E(f
\mid\F_\pi(\mathbf{s})).
\end{equation}

Consider two nonnegative random variables $Z_1$ and $Z_2$.
Then we define
%
\begin{equation}
Z_1 \unlhd Z_2
\quad \mbox{to mean}\quad
\P\bigl\{ \mathbf{1}_{\{Z_1>0\}}\le\mathbf{1}_{\{Z_2>0\}} \bigr\}=1,
\end{equation}
and $Z_1 \unrhd Z_2$ to mean $Z_2\unlhd Z_1$.
We also write $Z_1\asymp Z_2$ when $Z_1\unlhd Z_2$
and $Z_1\unrhd Z_2$. That is,
%
\begin{equation}
Z_1 \asymp Z_2 \quad \mbox{if and only if}\quad
\P\bigl\{ \mathbf{1}_{\{Z_1>0\}}=\mathbf{1}_{\{Z_2>0\}}
\bigr\}=1.
\end{equation}

The following generalizes Theorem 1.1 of Khoshnevisan and
Shi \cite{KS}. See also
Dalang and Nualart \cite{DN}, Theorem 31. This is the main
technical contribution of the present paper. We use the term
{\em upright box} for a Cartesian product $\Theta:=
\prod_{j=1}^N [a_j,b_j]$ of intervals, where $a_j < b_j$,
for $j=1,\dots,N$.

\begin{theorem}\label{thCPT}
Choose and fix an upright box $\Theta:=\prod_{j=1}^N[a_j,b_j]$
in $(0,\infty)^N$. For any partial order $\pi\subseteq\{1,
\ldots,N\}$, choose and fix some
vector $\mathbf{s}\in(0,\infty)^N\setminus\Theta$
such that $\mathbf{s} \prec_\pi\mathbf{t}$ for every $\mathbf
{t}\in\Theta$.
Then for all $\F_\pi(\mathbf{s})$-measurable bounded
random sets $A$,
%
\begin{equation}
\P^\pi_{\mathbf{s}} \{ B(\mathbf{u})\in A\mbox{ for some }
\mathbf{u}\in\Theta \} \asymp
\operatorname{Cap}_{d-2N}(A).
\end{equation}
\end{theorem}

We conclude this section with a technical result on
``potential theory of random sets.'' It should be
``obvious'' and/or well known. But we know of neither
transparent proofs nor explicit references. Therefore,
we supply a~proof.\looseness=-1\vadjust{\goodbreak}

\begin{lemma}\label{lemRPT}
Let $\mathcal{G}$ denote a~sub-$\sigma$-algebra
on the underlying probability space. Then for all
random $\mathcal{G}$-measurable closed sets $A\subseteq\R^d$
and all nonrandom $\beta\in\R$,
%
\begin{equation}\label{Eqcbeta}
\operatorname{Cap}_\beta(A)\asymp
[\inf\mathrm{I}_\beta(\theta) ]^{-1},
\end{equation}
where the infimum is taken over all random
$\mathcal{G}$-measurable probability measures $\theta$
that are compactly supported in $A$. In addition, there
is a $\mathcal{G}$-measurable random probability measure
$\mu$ such that $\operatorname{Cap}_\beta(A) \asymp1/ \mathrm{I}_\beta(\mu)$.
\end{lemma}

\begin{pf}
Let $c_\beta(A)$ denote the right-hand side of (\ref{Eqcbeta}).
Evidently,\break $\operatorname{Cap}_\beta(A)\ge c_\beta(A)$ almost surely,
and hence $\operatorname{Cap}_\beta(A) \unrhd c_\beta(A)$. It remains to
prove that $\operatorname{Cap}_\beta(A) \unlhd c_\beta(A)$. With this in
mind, we may---and will---assume without loss of generality that
$\operatorname{Cap}_\beta(A)>0$ with positive probability.
In particular, by \eqref{eqtaylor}, this implies that $\beta< d$.

Let $X_1,\ldots,X_M$ denote $M$ independent isotropic
stable processes in $\R^d$ that are
independent of the set $A$, and have a common stability index
$\alpha\in(0,2]$. Notice that we can always choose the integer
$M\ge1$ and the real number $\alpha$ such that
%
\begin{equation}\label{eqASP}
d-\alpha M = \beta.
\end{equation}
Thus, we choose and fix $(M,\alpha)$.

Define $\mathbf{X}$ to be the \emph{additive
stable process} defined by
%
\begin{equation}
\mathbf{X}(\mathbf{t}):= X_1(t_1)+\cdots+X_M(t_M)
\qquad \mbox{for all }\mathbf{t}\in\R^M_+,
\end{equation}
where we write $\mathbf{t}=(t_1,\ldots,t_M)$.
Theorem 4.1.1 of Khoshnevisan~\cite{Khbook}, page~423, tells us that for all nonrandom
\emph{compact} sets
$E\subseteq\R^d$,
\begin{eqnarray}
\P\{\mathbf{X}([1,2]^M)\cap E\neq\varnothing\}>0
&&\hspace*{-4pt}\quad \Longleftrightarrow\quad
\operatorname{Cap}_{d-\alpha M}(E)>0\nonumber\\[-8pt]\\[-8pt]
&&\hspace*{-4pt}\quad \Longleftrightarrow\quad
\operatorname{Cap}_\beta(E)>0;\nonumber
\end{eqnarray}
see \eqref{eqASP} for the final assertion.
The proof of that theorem (loc. cit.) tells us more.
Namely, that whenever $\operatorname{Cap}_\beta(E)>0$, there exists
a random variable~$\mathbf{T}$, with values in $[1,2]^M\cup
\{\infty\}$, which has the following properties:
\begin{itemize}[--]
\item[--] $\mathbf{T}\neq\infty$ if and only if $\mathbf{X}([1,2]^M)
\cap E\neq\varnothing$;
\item[--] $\mathbf{X}(\mathbf{T})\in E$ almost surely
on $\{\mathbf{T}\neq\infty\}$;
\item[--] $\mu(\bullet):=\P(\mathbf{X}(\mathbf{T})\in\bullet|
\mathbf{T}\neq\infty)$
is in $\mathcal{P}(E)$
and $\mathrm{I}_\beta(\mu)<\infty$.
\end{itemize}
In fact, $\mathbf{T}$ can be defined on $\{\mathbf{X}([1,2]^M)\cap
E\neq\varnothing\}$ as follows:
First define $T_1$ to be the smallest $s_1\in[1,2]$ such that
there exist $s_2,\ldots,s_M\in[1,2]$ that satisfy
$\mathbf{X}(s_1,\ldots,s_M)\in E$. Then, having defined $T_1,\ldots,T_j$
for $j\in\{1,\ldots,M-2\}$,
we define $T_{j+1}$ to be the smallest $s_{j+1}\in[1,2]$
such that there exist $s_{j+2},\ldots,s_M\in[1,2]$
that collectively satisfy
%
\begin{equation}
\mathbf{X}(T_1,\ldots,T_j,s_{j+1},\ldots,s_M)\in E.
\end{equation}
Finally, we define $T_M$ to be the smallest $s_M\in[1,2]$
such that
%
\begin{equation}
\mathbf{X}(T_1,\ldots,T_{M-1},s_M)\in E.\vadjust{\goodbreak}
\end{equation}
This defines $\mathbf{T}:=(T_1,\ldots,T_M)$
on $\{\mathbf{X}([1,2]^M)\cap E\neq\varnothing\}$.
We also define
$\mathbf{T}:=\infty$ on $\{
\mathbf{X}([1,2]^M)\cap E\neq\varnothing\}$. Then
$\mathbf{T}$ has the desired properties.

To finish the proof, note that, since $\operatorname{Cap}_\beta(A)>0$
with positive probability,
we can find $n>0$ such that $\operatorname{Cap}_\beta(A_n)>0$
with positive probability, where $A_n:=A\cap
[-n,n]^d$ is (obviously) a random $\mathcal{G}$-measurable compact set.
Because $A_n$ is independent of $\mathbf{X}$, we may apply the preceding
with $E:=A_n$. The mentioned construction of the resulting (now-random)
probability measure $\mu$ (on $A_n$) makes it clear
that $\mu$ is $\mathcal{G}$-measurable, and
$\mathrm{I}_\beta(\mu)<\infty$ almost surely
on $\{\operatorname{Cap}_\beta(A_n)>0\}$.
The lemma follows readily from these observations.
\end{pf}

\section{Analysis of pinned sheets}\label{sec3}

For all $\mathbf{s}\in(0,\infty)^N$ and $\mathbf{t}\in\R^N_+$, define
%
\begin{equation}\label{eqbridge}
B_{\mathbf{s}}(\mathbf{t}) := B(\mathbf{t}) - \delta_{\mathbf
{s}}(\mathbf{t}) B(\mathbf{s}),
\end{equation}
where
%
\begin{equation}\label{eqdelta}
\delta_{\mathbf{s}}(\mathbf{t}) := \prod_{j=1}^N \biggl(\frac{
s_j\wedge t_j}{s_j} \biggr).
\end{equation}
One can think of the random field $B_{\mathbf{s}}$ as the
\emph{sheet pinned to be zero at $\mathbf{s}$}.
(Khoshnevisan and Xiao \cite{KX} called $B_{\mathbf{s}}$ a ``bridge.'')

It is not too difficult to see that
%
\begin{equation}
B_{\mathbf{s}}(\mathbf{t}) = B(\mathbf{t}) - \E[ B(\mathbf{t}) |
B(\mathbf{s}) ].
\end{equation}

Next we recall some of the fundamental features of the pinned
sheet $B_{\mathbf{s}}$.

\begin{lemma}[(Khoshnevisan and Xiao \cite{KX}, Lemmas 51 and 52)]\label{lemcond}
Choose and fix a partial order $\pi\subseteq\{1,\ldots,N\}$
and a time point $\mathbf{s}\in(0,\infty)^N$.
Then $\{B_{\mathbf{s}}(\mathbf{t})\}_{\mathbf{t}\succ_\pi\mathbf{s}}$
is independent of
$\F_\pi(\mathbf{s})$. Moreover, for every nonrandom upright box
$I\subset(0,\infty)^N$ and $\pi\subseteq\{1,\ldots,N\}$,
there exists a finite constant $c>1$
such that uniformly for all $\mathbf{s}, \mathbf{u},\mathbf{v}\in I$,
%
\begin{equation}
c^{-1}\|\mathbf{u}-\mathbf{v}\|\le
\operatorname{Var} \bigl( B_{\mathbf{s}}^1(\mathbf{u})-B_{\mathbf
{s}}^1(\mathbf{v}) \bigr)
\le c\|\mathbf{u}-\mathbf{v}\|,
\end{equation}
where $B_{\mathbf{s}}^1(\mathbf{t})$ denotes the first coordinate of
$B_{\mathbf{s}}(\mathbf{t})$ for all $\mathbf{t}\in\R^N_+$.
\end{lemma}

The next result is the uniform Lipschitz continuity property of
the $\delta$'s.

\begin{lemma}\label{lemdeltamodulus}
Choose and fix an upright box $\Theta:=\prod_{j=1}^N[a_j,b_j]$.
Then there exists a constant $c< \infty$---depending only on
$N$, $\min_j a_j$ and $\max_j b_j$---such that
%
\begin{equation}
| \delta_{\mathbf{s}}(\mathbf{u})-\delta_{\mathbf{s}}(\mathbf{v}) |
\le c \|\mathbf{u}-\mathbf{v}\|
\qquad \mbox{ for all }\mathbf{s},\mathbf{u},\mathbf{v}\in\Theta.
\end{equation}
\end{lemma}

\begin{pf}
Notice that $\delta_s(t)$ is the product of $N$ bounded
and Lipschitz continuous functions $f_j(t_j) =
1\wedge(t_j/s_j)$, and the Lipschitz constants of these
functions are all bounded by $1/\min_j a_j$. The lemma follows.\vadjust{\goodbreak}
\end{pf}

Next, we present a conditional maximal inequality which extends
the existing multiparameter-martingale inequalities of the
literature in several directions.

\begin{lemma}\label{lemcairoli}
For every $\pi\subseteq\{1,\ldots,N\}$,
$\mathbf{s}\in\R^N_+$, and bounded $\sigma(B)$-measurable random
variable $f$,
%
\begin{equation}
\E^\pi_{\mathbf{s}} (\sup | \E^\pi_{\mathbf{t}} f
|^2 )\le4^N \E^\pi_{\mathbf{s}} (
|f|^2 ) \qquad \mbox{almost surely }[\P],
\end{equation}
where the supremum is taken over all $\mathbf{t}\in\mathbf{Q}^N_+$
such that $\mathbf{t}\succ_\pi\mathbf{s}$.
\end{lemma}

It is possible to use Lemma \ref{lemLBunif} below
in order to remove the restriction that $\mathbf{t}$
lies in $\mathbf{Q}^N_+$.

\begin{pf*}{Proof of Lemma \ref{lemcairoli}}
First we recall Cairoli's inequality,
%
\begin{equation}\label{eqcairoli}
\E\Bigl( \sup_{\mathbf{t}\in\mathbf{Q}^N_+} |
\E^\pi_{\mathbf{t}}f |^2 \Bigr) \le4^N
\E(|f|^2 ).
\end{equation}
When $\pi=\{1,\ldots,N\}$, this was proved by
Cairoli and Walsh \cite{CW}.
The general case is due to Khoshnevisan and Shi
\cite{KS}, Corollary 32. The proof
of \eqref{eqcairoli} hinges on the following
projection property (``commutation''):
%
\begin{equation}\label{eqcommutation}
\P\{ \E^\pi_{\mathbf{u}}\E^\pi_{\mathbf{t}}f
=\E^\pi_{\mathbf{u}\mathop{\curlywedge} _\pi\mathbf{t}} f \}=1,
\end{equation}
where, we recall, $\mathbf{u}\curlywedge_\pi\mathbf{t}$
denotes the $N$-vectors whose $j$th coordinate is
$u_j\wedge t_j$ if $j\in\pi$ and
$u_j\vee t_j$ if $j\notin\pi$.
Now we may observe that if $\mathbf{s}\prec_\pi\mathbf{u},\mathbf{t}$,
then $\P$-almost surely,
%
\begin{equation}
\P^\pi_{\mathbf{s}} \{ \E^\pi_{\mathbf{u}}\E^\pi_{\mathbf{t}}f
=\E^\pi_{\mathbf{u}\mathop{\curlywedge} _\pi\mathbf{t}} f \}=1.
\end{equation}
Thus, we apply the same proof that led to \eqref{eqcairoli},
but use the regular conditional distribution $\P^\pi_{\mathbf{s}}$
in place of $\P$, to finish the proof.
\end{pf*}

Next we mention a simple aside on certain Wiener integrals.

\begin{lemma}\label{lemgaussian}
Choose and fix a nonrandom compactly-supported bounded Borel
function 
$h\dvtx\R^d \to\R^d$, and a partial
order $\pi\subseteq\{1,\ldots,N\}$.
Define
%
\begin{equation}
G(\mathbf{s}):= \int_{\mathbf{r}\prec_\pi\mathbf{s}}h(\mathbf
{r}) B(\d\mathbf{r}),
\end{equation}
where the stochastic integral is defined in the sense of
Wiener \cite{Wiener23,Wiener38}, and~$\mathbf{s}$ ranges over $\R
^N_+$. Then $G$
has a continuous modification that is also continuous
in $L^2(\P)$.
\end{lemma}

\begin{pf}
Define $\mathfrak{S}(\mathbf{s})$ to be the $\pi$-shadow
of $\mathbf{s}\in\R^N_+$
%
\begin{equation}
\mathfrak{S}(\mathbf{s}):= \{ \mathbf{r}\in\R^N_+\dvtx
\mathbf{r}
\prec_\pi\mathbf{s} \}.\vadjust{\goodbreak}
\end{equation}
Then for all $\mathbf{s},\mathbf{t}\in\R^N_+$,
\begin{eqnarray}
\E\bigl( | G(\mathbf{t})-G(\mathbf{s}) |^2 \bigr)
&=&\int_{\mathfrak{S}(\mathbf{s})\triangle\mathfrak{S}(\mathbf{t})}
| h(\mathbf{r}) |^2 \,\d\mathbf{r}\nonumber\\[-8pt]\\[-8pt]
&\le&\sup|h|^2\times
\operatorname{meas} \bigl(\operatorname{supp}h\cap
\bigl(\mathfrak{S}(\mathbf{s})\triangle\mathfrak{S}
(\mathbf{t})\bigr) \bigr),\nonumber\vadjust{\goodbreak}
\end{eqnarray}
where ``$\operatorname{supp}h$'' denotes the support of $h$,
and ``$\operatorname{meas}$'' stands for the standard
$N$-dimensional Lebesgue measure.
Consequently,
$\E( | G(\mathbf{t})-G(\mathbf{s}) |^2 )
\le\mathrm{const}\cdot|\mathbf{s}-\mathbf{t}|$,
where the constant depends only on $(N,h)$. Because $G$ is
a Gaussian random field, it follows that
%
\begin{equation}
\E \bigl( | G(\mathbf{t})-G(\mathbf{s}) |^{2p} \bigr)
\le\mathrm{const}\cdot|\mathbf{s}-\mathbf{t}|^p
\qquad \mbox{for all }p>0,
\end{equation}
and the implied constant depends only on $(N,h,p)$.
The lemma follows
from a suitable form of the Kolmogorov continuity
lemma; see, for example, the arguments in
\v{C}encov \cite{Centsov} or Proposition A.1 and Remark
A.2 of Dalang et al.
\cite{DKN}.
\end{pf}

\begin{lemma}\label{lemLBunif}
Choose and fix a partial order $\pi\subseteq\{1,\ldots,N\}$.
If $Z$ is $\sigma(B)$-measurable and $\E(Z^2)<\infty$, then
$\mathbf{s}\mapsto\E^\pi_{\mathbf{s}} Z$ has a continuous
modification.
\end{lemma}

\begin{pf}
In the special case that $\pi=\{1,\dots,N\}$,
this is Proposition 2.3 of Khoshnevisan and Shi
\cite{KS}. Now we adapt
the proof to the present setting.

Suppose $h\dvtx\R^d\to\R$ is compactly supported and infinitely
differentiable. Define $\mathfrak{B}(h):=\int h\,\d B$, and note that
%
\begin{equation}\label{eqas}
\E^\pi_{\mathbf{s}} \bigl(\mathrm{e}^{\mathfrak{B}(h)} \bigr) =
\exp\biggl\{ \int_{\mathbf{r}\prec_\pi\mathbf{s}}
h(\mathbf{r}) B(\d\mathbf{r}) + \frac12\int_{\mathbf{r}\not\prec
_\pi\mathbf{s}}
| h(\mathbf{r}) |^2\,\d\mathbf{r} \biggr\}.
\end{equation}
Thanks to Lemma \ref{lemgaussian},
$\mathbf{s}\mapsto\E^\pi_{\mathbf{s}}[\exp(\mathfrak{B}(h))]$ is
continuous
almost surely. We claim that we also have continuity in $L^2(\P)$.
Indeed, we observe that it suffices to prove that $\mathbf{s}\mapsto
\exp(J_h(\mathbf{s}))$ is continuous in $L^2(\P)$, where
%
\begin{equation}
J_h(\mathbf{s}):=\int_{\mathbf{r}\prec_\pi\mathbf{s}}
h(\mathbf{r}) B(\d\mathbf{r}).
\end{equation}
By the Wiener isometry,
$\E(\exp(4J_h(\mathbf{s})))\le\exp(8\int|h(\mathbf{r})|^2\,\d
\mathbf{r})<\infty$.
By splitting the integral over $\mathbf{r}\prec_\pi\mathbf{s}$ into an
integral over $\mathbf{r} \in\mathfrak{S}(\mathbf{s})\setminus
\mathfrak{S}
(\mathbf{t})$ and a~remainder term, a direct calculation of $\mathrm{E}([\exp
(J_h(\mathbf{s})) - \exp(J_h(\mathbf{t}))]^2)$ using this inequality
yields the stated $L^2(\P)$ convergence.

We now use the preceding observation,
together with an approximation argument, as follows:

Thanks to Lemma 1.1.2 of Nualart (\cite{Nualart}, page 5), and by
the Stone--Weierstrass theorem, for all integers $m>0$,
we can find nonrandom compactly-supported functions
$h_1,\ldots,h_{k_m}\in C^\infty(\R^d)$ and $z_1,\ldots,z_{k_m}\in
\R$
such that
%
\begin{equation}\label{eqdense}
\E( | Z_m-Z |^2 )<\mathrm{e}^{-m}
\qquad \mbox{where }
Z_m:=\sum_{j=1}^{k_m} z_j \mathrm{e}^{\mathfrak{B}(h_j)}.
\end{equation}
Because conditional expectations are contractions on $L^2(\P)$,
it follows that
%
\begin{equation}\label{eqas1}
\E( | \E^\pi_{\mathbf{s}}Z-\E^\pi_{\mathbf{t}}Z |^2 )
\le9 \bigl(2\mathrm{e}^{-m}+\E( | \E^\pi_{\mathbf{s}}
Z_m - \E^\pi_{\mathbf{t}}
Z_m |^2 ) \bigr),
\end{equation}
and hence $\mathbf{s}\mapsto\E^\pi_{\mathbf{s}}Z$ is continuous in
$L^2(\P)$, therefore continuous in probability.

Thanks to \eqref{eqdense} and
Cairoli's maximal inequality \eqref{eqcairoli},
%
\begin{eqnarray}
\E\Bigl( \sup_{\mathbf{s}\in\mathbf{Q}^N_+}
| \E^\pi_{\mathbf{s}} Z_m
- \E^\pi_{\mathbf{s}}Z |^2 \Bigr)
&\le&4^N\sup_{\mathbf{s}\in\R^N_+}\E(
| \E^\pi_{\mathbf{s}}Z_m
- \E^\pi_{\mathbf{s}}Z |^2 )\nonumber \\
&\le&4^N \E(
|Z_m-Z |^2 )\\
&<& 4^N \mathrm{e}^{-m}.\nonumber
\end{eqnarray}
By the Borel--Cantelli lemma,
%
\begin{equation}
\lim_{m\to\infty}\sup_{\mathbf{s}\in\mathbf{Q}^N_+}
| \E^\pi_{\mathbf{s}} Z_m
- \E^\pi_{\mathbf{s}}Z | =0 \qquad \mbox{almost surely }[\P].
\end{equation}
Therefore, a.s. $[\P]$, the continuous random field $s \mapsto
\E^\pi_{\mathbf{s}} Z_m$ converges uniformly on $\mathbf{Q}^N_+$ to
$s \mapsto\E^\pi_{\mathbf{s}}Z$. Therefore, $s \mapsto\E^\pi
_{\mathbf{s}}Z$
is uniformly continuous on~$\mathbf{Q}^N_+$, and so it has a continuous
extension to $\R^N_+$. Since $s \mapsto\E^\pi_{\mathbf{s}}Z$ is continuous
in probability by \eqref{eqas1}, this extension defines a continuous
modification of~$\mathbf{s}\mapsto\E^\pi_{\mathbf{s}}Z$.
\end{pf}

Henceforth, we always choose a continuous modification of
$\E^\pi_{\mathbf{s}} Z$ when $Z$ is square-integrable. With this
convention in mind, we immediately obtain the following consequence
of Lemmas \ref{lemcairoli} and \ref{lemLBunif}.

\begin{lemma}\label{lemcairoli1}
For every bounded $\sigma(B)$-measurable random variable $f$, there
exists a
$\P$-null event off which the following holds:
For every $\pi\subseteq\{1,\ldots,N\}$,
$\mathbf{s}\in\R^N_+$,
%
\begin{equation}
\E^\pi_{\mathbf{s}} \Bigl(\sup_{\mathbf{t}\mathop{\succ}_\pi%
\mathbf{s}} | \E^\pi_{\mathbf{t}} f
|^2 \Bigr)\le4^N \E^\pi_{\mathbf{s}} (
|f|^2 ).
\end{equation}
\end{lemma}

For all $\sigma>0$, $\mathbf{t}\in\R^N$ and $z\in\R^d$ define
%
\begin{equation}\label{eqgamma}
\Gamma_{\sigma}(\mathbf{t} ; z) :=
\frac{1}{(2\pi\sigma^2)^{d/2} \|\mathbf{t}\|^{d/2}}
\exp\biggl( - \frac{\|z\|^2}{2\sigma^2\|\mathbf{t}\|} \biggr).
\end{equation}
Variants of the next result are well known.
We supply a detailed proof because we will need
to have good control over the constants involved.

\begin{lemma}\label{lemeasy}
Let $\Theta:=\prod_{j=1}^N[a_j,b_j]$ denote an upright box in
$(0,\infty)^N$,
and choose and fix positive constants
$\tau_1<\tau_2$ and $M>0$. Then there exists a finite constant\vadjust{\goodbreak}
$c>1$---depending only on $d$, $N$, $M$
$\tau_1$, $\tau_2$, $\min_j a_j$ and $\max_j b_j$---such that
for all $\sigma\in[\tau_1,\tau_2]$ and
$z\in[-M,M]^d$,
%
\begin{equation}\label{Eqkerbounds}
c^{-1}\kappa_{d-2N}(z)\le
\int_{\Theta- \Theta} \Gamma_\sigma(\mathbf{t} ;z)\, \d\mathbf{t}
\le c\kappa_{d-2N}(z).
\end{equation}
\end{lemma}

We recall that $\Theta-\Theta$ denotes the collection
of all points of the form $\mathbf{t}-\mathbf{s}$, where $\mathbf{s}$
and $\mathbf{t}$ range over $\Theta$. Moreover, the proof below shows
that the upper bound in (\ref{Eqkerbounds}) holds for all
$z\in\R^d$.

\begin{pf*}{Proof of Lemma \ref{lemeasy}}
Let $D(\rho)$ denote the centered ball in $\R^d$ whose
radius is $\rho>0$. Then we can integrate in polar
coordinates to deduce that
\begin{eqnarray}
\int_{D(\rho)} \Gamma_\sigma(\mathbf{t} ;z) \,\d\mathbf{t} &=&
\mathrm{const}\cdot\int_0^\rho r^{N-1-(d/2)} \exp\biggl(
-\frac{\|z\|^2}{2\sigma^2 r} \biggr)\, \d r\nonumber\\[-8pt]\\[-8pt]
&=& \frac{\mathrm{const}}{\|z\|^{d-2N}}\cdot\int_0^{2\sigma^2\rho/\|
z\|^2}
s^{N-1-(d/2)} \mathrm{e}^{-1/s}\, \d s,\nonumber
\end{eqnarray}
where the implied constants depend only on the parameters $\sigma$,
$N$ and $d$. This proves the result
in the case where $\Theta- \Theta$ is a centered ball, since
we can consider separately the cases $d<2N$, $d=2N$ and $d>2N$
directly; see the proof of Lemma~3.4 of
Khoshnevisan and Shi \cite{KS},
for instance.

The general case follows from the preceding
spherical case, because we can find $\rho_1$ and $\rho_2$
such that
$D(\rho_1)\subseteq\Theta- \Theta\subseteq
D(\rho_2)$, whence it follows that
$\int_{D(\rho_1)}\Gamma_\sigma(\mathbf{t} ;z) \,\d\mathbf{t}
\le\int_{\Theta- \Theta}\Gamma_\sigma(\mathbf{t} ;z) \,\d\mathbf{t}
\le\int_{D(\rho_2)}\Gamma_\sigma(\mathbf{t} ;z)\,\d\mathbf{t}$.
\end{pf*}

Now we proceed with a series of ``conditional energy estimates'' for
``continuous additive functionals'' of the sheet. First is a
lower bound.

\begin{lemma}\label{lemE1}
Choose and fix $\pi\subseteq\{1,\ldots,N\}$, $\eta>0$,
$\mathbf{s}\in(0,\infty)^N$, and
a~nonrandom upright box $\Theta:=\prod_{j=1}^N [a_j,b_j]$
in $(0,\infty)^N$ such that $\mathbf{s}\succ_\pi\mathbf{t}$
and $\eta\le|\mathbf{s}-\mathbf{t}|_\infty\le\eta^{-1}$ for
every $\mathbf{t}\in\Theta$.
Then there exists a constant $c>1$---depending only
on $d$, $N$, $\eta$, $\min_j a_j$ and $\max_j b_j$---such that
for all $\F_\pi(\mathbf{s})$-measurable
random probability density functions $f$ on $\R^d$,
%
\begin{equation}\label{eqE1}
\E^\pi_{\mathbf{s}} \biggl(\int_\Theta f(B(\mathbf{u}))\,\d\mathbf{u} \biggr)
\ge c^{-1} \mathrm{e}^{-c\|B(\mathbf{s})\|^2}\cdot\int_{\R^d} f(z)
\mathrm{e}^{-c\|z\|^2}\,\d z,
\end{equation}
almost surely $[\P]$.
\end{lemma}

\begin{pf}
Thanks to Lemma \ref{lemcond}, we can write
\begin{eqnarray}\label{eqid1}
\E^\pi_{\mathbf{s}}
\biggl(\int_\Theta f(B(\mathbf{u}))\,\d\mathbf{u} \biggr) &=& \E^\pi_{\mathbf
{s}} \biggl(
\int_\Theta f \bigl(
B_{\mathbf{s}} (\mathbf{u})+ \delta_{\mathbf{c}}(\mathbf
{u})B(\mathbf{s})
\bigr)\, \d\mathbf{u} \biggr)\nonumber\\[-8pt]\\[-8pt]
&=&\int_\Theta\d\mathbf{u}\int_{\R^d}\,\d z\,
f(z) g_{\mathbf{u}}\bigl(z-\delta_{\mathbf{s}}(\mathbf{u})B(\mathbf{s})\bigr),\nonumber
\end{eqnarray}
where $g_{\mathbf{u}}$ denotes the probability density function of
$B_{\mathbf{s}}(\mathbf{u})$.\vadjust{\goodbreak}

According to Lemma \ref{lemcond}, the coordinatewise variance
of $B_{\mathbf{s}}(\mathbf{u})$ is bounded above and
below by constant multiples of $\|\mathbf{u}- \mathbf{s}\|$.
As a result, $g_{\mathbf{u}}(z-\delta_{\mathbf{s}}(\mathbf
{u})B(\mathbf{s}))$ is bounded
below by an absolute constant multiplied by
\begin{eqnarray}
&&\frac{1}{\|\mathbf{u}- \mathbf{s}\|^{d/2}}
\exp\biggl( -\mathrm{const} \frac{
\|z-\delta_{\mathbf{s}}(\mathbf{u})B(\mathbf{s})\|^2} {\|\mathbf
{u}- \mathbf{s}\|} \biggr)\nonumber\\[-8pt]\\[-8pt]
&&\qquad
\ge\eta^{d/2}\exp\biggl(-\mathrm{const}\cdot
\frac{\|z-\delta_{\mathbf{s}}(\mathbf{u})B(\mathbf{s})\|^2}{\eta} \biggr).\nonumber
\end{eqnarray}
Thus, the inequality
%
\begin{equation}
\|z-\delta_{\mathbf{s}}(\mathbf{u})B(\mathbf{s})\|^2\le
2\|z\|^2+2\|B(\mathbf{s})\|^2,
\end{equation}
valid because $0\le\delta_{\mathbf{s}}(\mathbf{u})\le1$,
proves that
%
\begin{equation}\label{eqdensityLB}
g_{\mathbf{u}}\bigl(z-\delta_{\mathbf{s}}(\mathbf{u})B(\mathbf{s})\bigr)\ge
c_1 \exp\bigl( -c_2 \{
\|z\|^2 +\|B(\mathbf{s})\|^2 \} \bigr),
\end{equation}
where $c_1$ and $c_2$ are positive
and finite constants that depend only on
$\pi$, $d$, $N$, $M$, $\eta$ and $\max_j b_j$. Let $c_1(\pi)$
and $c_2(\pi)$ denote the same constants, but written
as such to exhibit their dependence on the partial order
$\pi$. Apply the
preceding for all partial orders $\pi$, and
let $c_1$ and $c_2$ denote, respectively, the minimum and
maximum of $c_1(\pi)$ and $c_2(\pi)$ as $\pi$
ranges over the various subsets of $\{1,\ldots,N\}$.
In this way, the preceding display holds
without any dependencies on the partial order $\pi$.
It is now clear that (\ref{eqE1}) follows from~(\ref{eqid1})
and (\ref{eqdensityLB}).
\end{pf}

Next we present a delicate joint-density estimate
for the pinned sheets. This estimate will be
used subsequently to describe a conditional second-moment
bound that complements
the conditional
first-moment bound of Lemma \ref{lemE1}.

\begin{lemma}\label{lemjointpdf}
Choose and fix an upright box $\Theta:=\prod_{j=1}^N
[a_j,b_j]$ in $(0,\infty)^N$, a partial
order $\pi\subseteq\{1,\ldots,N\}$ and $\mathbf{s}\in\R^N_+$
and $\eta>0$ such that:
\begin{longlist}[(ii)]
\item[(i)] $\mathbf{s}\prec_\pi\mathbf{t}$ for all $\mathbf{t}\in
\Theta$;
\item[(ii)] $\eta\le|\mathbf{s}-\mathbf{t}|_\infty\le\eta^{-1}$
for all $\mathbf{t}\in\Theta$.
\end{longlist}
Then there exists a finite constant $c>1$---depending
only on $d$, $N$, $\eta$, $\min_j a_j$ and $\max_j b_j$---such that
for all
$x,y\in\R^d$ and $\mathbf{u},\mathbf{v}\in\Theta$,
%
\begin{equation}
p_{\mathbf{s};\mathbf{u},\mathbf{v}}(x,y) \le
c \Gamma_c(\mathbf{u}-\mathbf{v} ;x-y),
\end{equation}
where $p_{\mathbf{s};\mathbf{u},\mathbf{v}}(x,y)$ denotes
the probability density function of $(B_{\mathbf{s}}(\mathbf{u})
,B_{\mathbf{s}}(\mathbf{v}))$.
\end{lemma}

\begin{pf}
The proof is carried out in three steps. We are only going to consider
the case where
$u \neq u \curlywedge_\pi v \neq v$; indeed, the other
cases are simpler and are left to the reader.

\textit{Step 1}. First consider the case that $\pi=\{1,\ldots,N\}$.
In this particular case, we, respectively, write ``$\prec$,''
``$\succ$'' and ``$\curlywedge$'' in place of ``$\prec_\pi$,''
``$\succ_\pi$''
and ``$\curlywedge_\pi$.''\vadjust{\goodbreak}

Note that $\mathbf{r}\prec\mathbf{p}$ if and only if
$r_i\le p_i$ for all $i=1,\ldots,N$. Furthermore,
%
\begin{equation}
B_{\mathbf{s}}(\mathbf{r}) = B(\mathbf{r})-B(\mathbf{s})
\qquad \mbox{for all }\mathbf{r}\in\Theta.
\end{equation}

Because the joint probability-density function of $(B_{\mathbf
{s}}(\mathbf{u})
,B_{\mathbf{s}}(\mathbf{v}))$ is unaltered if we modify
the Brownian sheet, we choose to work with a particularly
useful construction of the Brownian sheet. Namely,
let $\mathfrak{W}$ denote $d$-dimensional white noise
on $\R^N_+$, and consider the Brownian sheet
%
\begin{equation}
B(\mathbf{t}):=\mathfrak{W}([\mathbf{0},\mathbf{t}]),
\qquad \mbox{where }[\mathbf{0},\mathbf{t}]:=\prod_{j=1}^N[0,t_j].
\end{equation}
This construction might not yield a continuous random
function $B$, but that is not germane to the discussion.

For the construction cited here,
%
\begin{equation}
B_{\mathbf{s}}(\mathbf{r}) = \mathfrak{W} ( [\mathbf{0},\mathbf{r}]
\setminus[\mathbf{0},\mathbf{s}] )
\qquad \mbox{for all }\mathbf{r}\in\Theta.
\end{equation}

For all bounded $C^\infty$ functions $\phi\dvtx(\R^d)^2\to\R_+$
and $\mathbf{u},\mathbf{v}\in\Theta$,
%
\begin{equation}
\E[\phi( B_{\mathbf{s}}(\mathbf{u}) ,
B_{\mathbf{s}}(\mathbf{v}) ) ] =\iiint
\phi(x+y,x+z)g_{\mathbf{u}\curlywedge\mathbf{v}}(x) F(y) G(z)\,
\d x\,\d y\,\d z,\hspace*{-40pt}
\end{equation}
where $g_{\mathbf{u}\curlywedge\mathbf{v}}$ denotes the probability density
function of $B_{\mathbf{s}}(\mathbf{u}\curlywedge\mathbf{v})
=\mathfrak{W}([\mathbf{0},\mathbf{u}\curlywedge\mathbf
{v}]\setminus[\mathbf{0},\mathbf{s}])$
as before, $F$ the probability
density function of $\mathfrak{W}([\mathbf{0},\mathbf{u}]
\setminus[\mathbf{0},\mathbf{u}\curlywedge\mathbf{v}])$ and $G$ the
probability density function of
$\mathfrak{W}([\mathbf{0},\mathbf{v}]
\setminus[\mathbf{0},\mathbf{u}\curlywedge\mathbf{v}])$. The
integrals are each taken over
$\R^d$. The $N$-dimensional volume of
$[\mathbf{0},\mathbf{u}\curlywedge\mathbf{v}]\setminus[\mathbf
{0},\mathbf{s}]$ is
at least $\eta(\min_j a_j)^{N-1}$. Therefore,
$g_{\mathbf{u}\curlywedge\mathbf{v}}$ is bounded above by a constant
$c_3$ that depends only
on $d$, $N$, $\eta$ and $\min_j a_j$. And hence,
\begin{eqnarray}\qquad
\E[\phi( B_{\mathbf{s}}(\mathbf{u}) ,
B_{\mathbf{s}}(\mathbf{v}) ) ] &\le& c_3 \iiint
\phi(x+y,x+z) F(y) G(z)\,
\d x\,\d y\,\d z\nonumber\\[-8pt]\\[-8pt]
&=&c_3\iint\phi(x,y) (F*G)(y-x)\,\d x\,\d y.\nonumber
\end{eqnarray}
But $F*G$ is the probability density function of
$\mathfrak{W}([\mathbf{0},\mathbf{u}]\,
\triangle\,[\mathbf{0},\mathbf{v}])$, and the $N$-dimensional
volume of $[\mathbf{0},\mathbf{u}]\,
\triangle\,[\mathbf{0},\mathbf{v}]$ is at least
%
\begin{equation}\qquad
\Bigl( \min_{1\le j\le N} a_j \Bigr)^{N-1}
\sum_{k=1}^N|u_k-v_k|
\ge
\frac{1}{N^{1/2}} \Bigl( \min_{1\le j\le N} a_j
\Bigr)^{N-1}\|\mathbf{u}-\mathbf{v}\|.
\end{equation}
In addition, one can derive an upper bound---using only constants that
depend on $\min_j a_j$, $\max_j b_j$ and $N$---similarly.
Therefore, there exists a finite constant $c>1$---depending only
on $d$, $N$, $\eta$ and $\min_j a_j$---such that the following
occurs pointwise:
%
\begin{equation}
(F*G)(y-x)\le c \Gamma_{c}(\mathbf{u}-\mathbf{v} ;y-x).
\end{equation}
This proves the lemma in the case where
$\pi=\{1,\ldots,N\}$.\vadjust{\goodbreak}

\textit{Step 2}. The argument of Step 1 yields in fact a slightly
stronger result, which we state next as
the following (slightly)
\textbf{Enhanced Version}: \emph{Choose and fix
two positive constants $\nu_1<\nu_2$. Under the conditions of
Step 1, there exists a constant $\rho$---depending only
on $d$, $N$, $\eta$, $\min_j a_j$, $\max_j b_j$,
$\nu_1$ and $\nu_2$, such that for all
$\mathbf{u},\mathbf{v}\in\Theta$ and all $\alpha,\beta\in[\nu
_1,\nu_2]$,
the joint probability density function
of $(\alpha B_{\mathbf{s}}(\mathbf{u}),
\beta B_{\mathbf{s}}(\mathbf{v}))$---at $(x,y)$---is bounded
above by $\rho\Gamma_\rho(\mathbf{u}-\mathbf{v};\allowbreak x-y)$.}\looseness=1

The proof of the enhanced version is
the same as the case we expanded on above ($\nu_1=\nu_2=1$).
However, a few modifications need to be made: $\phi(x+y,x +z)$
is replaced by $\phi(\alpha x+y,\beta x +z)$;
$F$ is replaced by the probability density function of
$\alpha\mathfrak{W}([\mathbf{0},\mathbf{u}]\setminus[\mathbf{0},
\mathbf{u}\curlywedge\mathbf{v}])$; $G$ by the probability density
function of
$\beta\mathfrak{W}([\mathbf{0},\mathbf{v}]\setminus[\mathbf{0},
\mathbf{u}\curlywedge\mathbf{v}])$; and $F*G$ is now the probability density
function of a centered Gaussian
vector with i.i.d. coordinates, the variance of
each of which is at least
\begin{eqnarray}
\qquad && \Bigl( \min_{1\le j\le N} a_j \Bigr)^{N-1}\alpha^2
\sum_{k=1}^N(u_k-v_k)^+ + \Bigl(
\min_{1\le j\le N} a_j \Bigr)^{N-1}\beta^2
\sum_{k=1}^N(u_k-v_k)^- \nonumber\\[-8pt]\\[-8pt]
\qquad &&\qquad  \ge(\alpha\wedge\beta)^2
\Bigl( \min_{1\le j\le N}a_j \Bigr)^{N-1}
\sum_{k=1}^N|u_k-v_k|.\nonumber
\end{eqnarray}
The remainder of the proof goes through
without incurring major changes.

\textit{Step 3}. If $\pi=\varnothing$, then the lemma follows
from Step 1 and symmetry. Therefore,
it remains to consider the case where
$\pi$ and $\{1,\ldots,N\}\setminus\pi$ are both nonvoid. We follow
Khoshnevisan and Xiao \cite{KX}, proof of Proposition 31,
and define a map $\mathcal{I}\dvtx(0,\infty)^N
\to(0,\infty)^N$ with coordinate functions
$\mathcal{I}_1,\ldots,\mathcal{I}_N$ as follows:
For all $k=1,\ldots,N$,
%
\begin{equation}\label{eqcalI}
\mathcal{I}_k(\mathbf{t}) :=
\cases{
t_k,&\quad if $k\in\pi$,\cr
1/t_k,&\quad if $k\notin\pi$.
}
\end{equation}

Consider any two points $\mathbf{u},\mathbf{v}\in\Theta$.
We may note that:
\begin{itemize}[(iii)]
\item[(i)] $\mathcal{I}(\Theta)$ is an
upright box that contains $\mathcal{I}(\mathbf{u})$
and $\mathcal{I}(\mathbf{v})$;
\item[(ii)] $\mathcal{I}(\mathbf{s})\prec\mathcal{I}(\mathbf{t})$
for all $\mathbf{t}\in\Theta$ (nota bene: the partial order!);
\item[(iii)] $|\mathcal{I}(\mathbf{s})-\mathcal{I}(\mathbf
{t})|_\infty$
is bounded below by a positive constant $\eta'$, uniformly
for all $\mathbf{t}\in\Theta$. Moreover, $\eta'$ depends only
on $N$, $\eta$, $\min_j a_j$ and $\max_j b_j$.
\end{itemize}
Define
%
\begin{equation}
W(\mathbf{t}) := \biggl(\prod_{j\notin\pi} t_j \biggr) \cdot
B(\mathcal{I}(\mathbf{t}))
\qquad \mbox{for all }\mathbf{t}\in(0,\infty)^N.
\end{equation}
Then, according to Khoshnevisan and Xiao (loco citato),
$W$ is a Brownian sheet. Thus, we have also the corresponding
pinned sheet
%
\begin{equation}
W_{\mathbf{s}}(\mathbf{t}) = W(\mathbf{t}) - \delta_{\mathbf
{s}}(\mathbf{t})
W(\mathbf{s}) \qquad \mbox{for all }\mathbf{t}\in(0,\infty)^N.
\end{equation}
It is the case that
%
\begin{equation}\qquad
W_{\mathbf{s}}(\mathbf{t}) = \biggl(\prod_{j\notin\pi} t_j \biggr) \cdot
[ B(\mathcal{I}(\mathbf{t}))-
B(\mathcal{I}(\mathbf{s})) ]
\qquad \mbox{for all }\mathbf{t}\in(0,\infty)^N.
\end{equation}
The derivation of this identity requires
only a little algebra, which we skip. Thus,
property (ii) above implies the following remarkable
identity:
%
\begin{equation}
W_{\mathbf{s}}(\mathbf{t}) = \biggl(\prod_{j\notin\pi} t_j \biggr) \cdot
B_{\mathcal{I}(\mathbf{s})}(\mathcal{I}(\mathbf{t}))
\qquad \mbox{for all }\mathbf{t}\in\Theta.
\end{equation}

As a result of items (i)--(iii), and thanks to Step 1,
the joint probability density function---at $(x,y)$---of the
random vector $(B_{\mathcal{I}(\mathbf{s})}
(\mathcal{I}(\mathbf{u})),B_{\mathcal{I}(\mathbf{s})}(\mathcal
{I}(\mathbf{v})))$
is bounded above by $c_4\Gamma_{c_4}
(\mathcal{I}(\mathbf{u})-\mathcal{I}(\mathbf{v}) ;x-y)$,
where $c_4$ depends only on~$d$, $N$, $\eta$,
$\min_ja_j$ and $\max_j b_j$. Elementary considerations
show that $\|\mathcal{I}(\mathbf{u})-\mathcal{I}(\mathbf{v})\|$
is bounded above and below by constant multiples of
$\|\mathbf{u}-\mathbf{v}\|$, where the constants have the
same parameter dependencies as $c_4$.
These discussions together imply that the joint probability density
function---at $(x,y)$---of the random vector
$(B_{\mathcal{I}(\mathbf{s})}
(\mathcal{I}(\mathbf{u})),B_{\mathcal{I}(\mathbf{s})}(\mathcal
{I}(\mathbf{v})))$ is bounded
above by $c_5\Gamma_{c_5}(\mathbf{u}-\mathbf{v} ;\allowbreak x-y)$,
where $c_5$ has the same parameter dependencies
as $c_4$. Set $\alpha=\prod_{j\notin\pi}u_j$ and
$\beta:=\prod_{j\notin\pi}v_j$, and note that
$\alpha$ and $\beta$ are bounded above and below by
constants that depend only on $\min_j a_j$ and $\max_j b_j$.
Also note that
%
\begin{equation}
\bigl(\alpha B_{\mathcal{I}(\mathbf{s})}
(\mathcal{I}(\mathbf{u})) , \beta B_{\mathcal{I}(\mathbf
{s})}(\mathcal{I}(\mathbf{v}))\bigr)
= (W_{\mathbf{s}} (\mathbf{u}),W_{\mathbf{s}}(\mathbf{v})),
\end{equation}
for all $\mathbf{u},\mathbf{v}\in\Theta$.
Thus, in accord with Step 2, the joint probability
density function---at $(x,y)$---of
$(W_{\mathbf{s}} (\mathbf{u}),W_{\mathbf{s}}(\mathbf{v}))$ is bounded
above by $c_6\Gamma_{c_6}(\mathbf{u}-\mathbf{v} ;\allowbreak x-y)$, where
$c_6$ has the same parameter dependencies as $c_4$.
Because $W_{\mathbf{s}}$ has the same finite-dimensional
distributions as $B_{\mathbf{s}}$, this proves the lemma.~%
\end{pf}

\begin{lemma}\label{lemE2}
Let $\Theta$, $\mathbf{s}$, $\pi$ and $\eta$ be as in
Lemma \ref{lemjointpdf}.
Then there exists a constant $c>1$---depending only
on $d$, $N$, $\eta$, $\min_j a_j$ and $\max_j b_j$---such that
for all $\F_\pi(\mathbf{s})$-measurable random probability
density functions $f$,
%
\begin{equation}\quad
\E^\pi_{\mathbf{s}} \biggl( \biggl|
\int_\Theta f(B(\mathbf{u}))\,\d\mathbf{u} \biggr|^2 \biggr)
\le c \mathrm{e}^{c\|B(\mathbf{s})\|^2} \cdot\mathrm{I}_{d-2N}(f)
\qquad \mbox{a.s. }[\P].
\end{equation}
\end{lemma}

\begin{pf}
Throughout this proof, we define
%
\begin{equation}
F := \E^\pi_{\mathbf{s}} \biggl( \biggl|\int_\Theta
f(B(\mathbf{u}))\,\d\mathbf{u} \biggr|^2 \biggr).
\end{equation}
A few lines of computation show that with probability one,
\begin{eqnarray}\qquad
F&=&\int_{\Theta}\,\d\mathbf{v}\int_{\Theta}\,\d\mathbf{u}\int_{\R
^d}\,\d x
\int_{\R^d}\,\d y\, f\bigl(x+\delta_{\mathbf{s}}(\mathbf{u})B(\mathbf{s})\bigr)\nonumber\\[-8pt]\\[-8pt]
&&\hphantom{\int_{\Theta}\,\d\mathbf{v}\int_{\Theta}\,\d\mathbf{u}\int_{\R
^d}\,\d x
\int_{\R^d}\,\d y\,}
{}\times f\bigl(y+\delta_{\mathbf{s}}(\mathbf{v})B(\mathbf{s})\bigr) p_{\mathbf
{s};\mathbf{u},\mathbf{v}}(x,y),\nonumber
\end{eqnarray}
where $p_{\mathbf{s};\mathbf{u},\mathbf{v}}(x,y)$ denotes the
probability density function of
$(B_{\mathbf{s}}(\mathbf{u}),B_{\mathbf{s}}(\mathbf{v}))$ at
$(x,y)\in(\R^d)^2$.
According to Lemma \ref{lemjointpdf}, we can find a
finite constant $c_7>1$ such that for all
$(x,y)\in(\R^d)^2$ and $\mathbf{u},\mathbf{v}\in\Theta$,
%
\begin{equation}
p_{\mathbf{s};\mathbf{u},\mathbf{v}}(x,y) \le
c_7 \Gamma_{c_7}(\mathbf{u}-\mathbf{v} ; x-y),
\end{equation}
where $\Gamma_c$ is the Gaussian density function defined by \eqref{eqgamma}.
Moreover, $c_7$ depends only on $d$, $M$, $N$,
$\eta$, $\min_j a_j$ and $\max_j b_j$. We change variables
to deduce that almost surely
%
\begin{equation}
F \le c_7\int_{\Theta}\d\mathbf{v}\int_{\Theta}\d\mathbf{u}\int
_{\R^d}\,\d x
\int_{\R^d}\,\d y\,f(x) f(y)
\Gamma_{c_7} ( \mathbf{u}-\mathbf{v} ; x-y-
Q ),\hspace*{-35pt}
\end{equation}
where
%
\begin{equation}
Q:=B(\mathbf{s})\{\delta_{\mathbf{s}}(\mathbf{u})-
\delta_{\mathbf{s}}(\mathbf{v})\}.
\end{equation}
Because $\|z\|^2\le2\|Q\|^2+2\|z-Q\|^2$,
%
\begin{equation}
\Gamma_{c_7} (\mathbf{t} ; z-Q) \le
\frac{1}{(2\pi c_7^2)^{d/2} \|\mathbf{t}\|^{d/2}}
\exp\biggl( - \frac{\|z\|^2}{c_7\|\mathbf{t}\|}
+\frac{\|Q\|^2}{c_7\|\mathbf{t}\|} \biggr).
\end{equation}
According to Lemma \ref{lemdeltamodulus}, there exists a constant $c_8$---with
the same parameter dependencies as $c_7$---such that
\begin{eqnarray}
\|Q\|^2&\le& c_8\|\mathbf{u}-\mathbf{v}\|^2 \cdot\|B(\mathbf{s})\|
^2\nonumber\\[-8pt]\\[-8pt]
&\le& c_8\max_{1\le j\le N} b_j\|\mathbf{u}-\mathbf{v}\|\cdot\|
B(\mathbf{s})\|^2,\nonumber
\end{eqnarray}
uniformly for all $\mathbf{u},\mathbf{v}\in\Theta$.
Therefore, we may apply the preceding display with
$\mathbf{t}:=\mathbf{u}-\mathbf{v}$ and $z:=x-y$ to find that
%
\begin{equation}
\Gamma_{c_7} (\mathbf{u}-\mathbf{v} ; x-y-Q) \le2^d\exp\biggl(
\frac{\|B(\mathbf{s})\|^2}{c_7c_8\max_{1\le j\le N} b_j} \biggr)\,{\cdot}\,
\Gamma_{c_9}(\mathbf{u}-\mathbf{v} ; x-y).\hspace*{-40pt}
\end{equation}
Again, $c_9$ is a positive and finite constant that has the
same parameter dependencies as $c_7$ and $c_8$.
Consequently, the following holds with probability one:
%
\begin{eqnarray}\label{eqwtbd}
&&\exp\biggl(-\frac{\|B(\mathbf{s})\|^2}{c_7c_8\max_{1\le j\le N} b_j} \biggr)\cdot
F\nonumber \\
&&\qquad \le2^dc_7\int_{\Theta}\,\d\mathbf{v}\int_{\Theta}\d
\mathbf{u}\int_{\R^d}
\d x\int_{\R^d}\,\d y\, f(x) f(y)
\Gamma_{c_9} ( \mathbf{u}-\mathbf{v} ; x-y
)\\
&&\qquad = 2^dc_7\operatorname{meas}(\Theta)\int_{\R^d}\d x\int
_{\R^d}\,\d y\, f(x) f(y) g(x-y),\nonumber
\end{eqnarray}
where $g(z) := \int_{\Theta- \Theta}
\Gamma_{c_9} ( \mathbf{u} ; z)\,\d\mathbf{u}$. Thanks to
Lemma \ref{lemeasy},
%
\begin{equation}
g(z)\le c_{10} \kappa_{d-2N}(z)
\end{equation}
for all $z\in\R^d$, where $c_{10}$ is a finite constant $>1$ that
depends only on $d$, $N$, $M$, $\eta$, $\min_j a_j$ and $\max_j
b_j$.\vadjust{\goodbreak}
The lemma follows.
\end{pf}

Next we introduce a generalization of Proposition 3.7 of
Khoshnevisan and Shi \cite{KS}.

\begin{lemma}\label{lemLB}
Choose and fix an upright box $\Theta:=
\prod_{j=1}^N[a_j,b_j]$ in $(0,\infty)^N$
and real numbers $\eta>0$ and $M > 0$. Then there exists
a constant $c_{11}>0$---depending only
on $d$, $N$, $\eta$, $M$, $\min_j a_j$ and $\max_j b_j$---such
that for all
$\pi\subseteq\{1,\ldots,N\}$, all $\mathbf{s}\in\Theta$ whose distance
to the boundary of $\Theta$ is at least~$\eta$, and every
$\F_\pi(\mathbf{s})$-measurable random probability density function
$f$ whose support is contained in $[-M, M]^N$,
\begin{eqnarray}\label{eqLB}
&&\E^\pi_{\mathbf{s}} \biggl(
\int_\Theta f(B(\mathbf{u}))\,\d\mathbf{u}
\biggr)\nonumber\\[-8pt]\\[-8pt]
&&\qquad \ge c_{11} \1_{\{\|B(s)\|\le M\}}
\int_{\R^d}\kappa_{d-2N}(z)f\bigl(z+B(\mathbf{s})\bigr)\,\d z,\nonumber
\end{eqnarray}
almost surely $[\P]$.
\end{lemma}

Even though both Lemmas \ref{lemE1} and \ref{lemLB} are
concerned with lower bounds for $\E^\pi_{\mathbf{s}} (
\int_\Theta f(B(\mathbf{u}))\,\d\mathbf{u} )$, there is a
fundamental difference between the two lemmas: in Lemma \ref{lemE1},
$\mathbf{s}$ is at least a fixed distance $\eta$ away from $\Theta$,
whereas Lemma \ref{lemLB} considers the case where $\mathbf{s}$
belongs to $\Theta$.

\begin{pf*}{Proof of Lemma \ref{lemLB}}
Throughout, we choose and fix an $\mathbf{s} \in\Theta$ and
a $\pi$ as per the statement of the lemma.

Consider
$\Upsilon:=\{\mathbf{u}\in\Theta\dvtx \mathbf{u}\succ_\pi\mathbf
{s}\}$,
which is easily seen to be an upright box.
Since $\Upsilon\subseteq\Theta$, it follows that $\P$-almost surely,
%
\begin{eqnarray}
\E^\pi_{\mathbf{s}} \biggl(\int_\Theta f(B(\mathbf{u}))\,\d\mathbf{u} \biggr)
&\ge&\E^\pi_{\mathbf{s}} \biggl(\int_{\Upsilon} f(B(\mathbf{u}))\,\d
\mathbf{u} \biggr)\nonumber \\
&=&\E^\pi_{\mathbf{s}} \biggl(\int_{\Upsilon} f \bigl( B_{\mathbf{s}}(\mathbf{u})
+\delta_{\mathbf{s}}(\mathbf{u})B(\mathbf{s}) \bigr)\,\d\mathbf{u} \biggr)\\
&=&\int_{\Upsilon}\d\mathbf{u}\int_{\R^d}\,\d z\, f(z)
g_{\mathbf{u}} \bigl( z-\delta_{\mathbf{s}}(\mathbf{u})B(\mathbf{s}) \bigr),\nonumber
\end{eqnarray}
where $g_{\mathbf{u}}$ denotes the probability
density function of $B_{\mathbf{s}}(\mathbf{u})$,
as before. We temporarily use the abbreviated notion
$\delta:=\delta_{\mathbf{s}}(\mathbf{u})$ and $y:=B(\mathbf{s})$.
Thanks to Lemma~\ref{lemcond}, for all $z\in\R^d$,
%
\begin{equation}
g_{\mathbf{u}}(z-\delta y) \ge\frac{c_{12}}{\|\mathbf{u}-\mathbf
{s}\|^{d/2}}
\exp\biggl( -\frac{\|z-\delta y\|^2}{c_{12}\|\mathbf{u}-\mathbf{s}\|} \biggr),
\end{equation}
where $c_{12}\in(0,1)$ depends only on $N$, $\min_j a_j$ and
$\max_j b_j$. But
%
\begin{equation}
\|z-\delta y\|^2\le
2\|z-y\|^2+2\|y\|^2(1-\delta)^2
\end{equation}
and
%
\begin{equation}
0\le1-\delta= \delta_{\mathbf{s}}(\mathbf{s})-\delta_{\mathbf
{s}}(\mathbf{u})
\le\mathrm{const}\cdot\|\mathbf{u}-\mathbf{s}\|,\vadjust{\goodbreak}
\end{equation}
for a constant that has
the same parameter dependencies as $c_{12}$. Consequently,
there exists $c_{13}\in(0,1)$---depending
only on $N$, $\min_j a_j$ and $\max_j b_j$---such that
%
\begin{equation}
g_{\mathbf{u}}(z-\delta y)\ge c_{13}
\mathrm{e}^{-\|B(\mathbf{s})\|^2/c_{13}}
\cdot\Gamma_{c_{13}}(\mathbf{u}-\mathbf{s} ;
z-y).
\end{equation}
Recall that $\delta=\delta_{\mathbf{s}}(\mathbf{u})$
and $y:=B(\mathbf{s})$; it follows from this discussion that
$\P$-almost surely,
\begin{eqnarray}
&&\E^\pi_{\mathbf{s}} \biggl(\int_\Theta f(B(\mathbf{u}))\,\d\mathbf{u}
\biggr)\nonumber\\[-8pt]\\[-8pt]
&&\qquad  \ge c_{14} \mathrm{e}^{-\|B(\mathbf{s})\|^2/c_{14}}
\int_{\R^d}\,\d z\, f\bigl(z+B(\mathbf{s})\bigr)
\biggl(\int_{\Upsilon-\mathbf{s}}\Gamma_{c_{14}}(\mathbf{u} ;z)\,\d
\mathbf{u}
\biggr).\nonumber
\end{eqnarray}
Because the distance between $\mathbf{s}$ and the boundary of
$\Theta$ is at least $\eta$, the upright box $\Upsilon-\mathbf{s}$ contains
$[0,\eta]^N$. Therefore, by symmetry,
\begin{eqnarray}
\qquad &&\E^\pi_{\mathbf{s}} \biggl(\int_\Theta f(B(\mathbf{u}))\,\d\mathbf{u}
\biggr)\nonumber\\[-8pt]\\[-8pt]
\qquad &&\qquad \ge\frac{c_{13}}{2^N} \mathrm{e}^{-\|B(\mathbf{s})\|^2/c_{13}}
\int_{\R^d}\,\d z\, f\bigl(z+B(\mathbf{s})\bigr)
\biggl(\int_{[-\eta,\eta]^N}\Gamma_{c_{13}}(\mathbf{u} ;z)\,\d\mathbf{u}
\biggr),\nonumber
\end{eqnarray}
almost surely $[\P]$. Since the support of $f$ is contained in $[-M, M]^N$,
Lem\-ma~\ref{lemeasy} finishes the proof.
\end{pf*}

\section{\texorpdfstring{Proof of Theorem \protect\ref{thCPT}}{Proof of Theorem 2.4}}\label{sec4}
We begin by making two simplifications:
\begin{itemize}
\item First, let us note that the upright box $\Theta$ is closed, and hence
there exists $\eta\in(0, 1)$ such that $\eta\le|\mathbf{s}-\mathbf
{t}|_\infty
\le\eta^{-1}$ for all $\mathbf{t}\in\Theta$. This $\eta$ is held
fixed throughout the proof.
\item Thanks to the capacitability theorem of Choquet, we
may consider only $\F_\pi(\mathbf{s})$-measurable
\emph{compact} random sets $A\subset[-M,M]^d$.
Without loss of generality, we may---and will---assume
that $M>1$ is fixed henceforth.
\end{itemize}

For every nonrandom $\epsilon\in(0,1)$,
we let $A^\epsilon$ denote the closed $\epsilon$-enlargement
of~$A$. Let $f$ denote a random
$\F_\pi(\mathbf{s})$-measurable density function
that is supported on $A^\epsilon$.
Because we assumed that $M$ is greater than one,
and $\epsilon$ is at most one,
$\|z\|^2\le\mathrm{const}\cdot M^2$
for all $z\in A^\epsilon$ and $\epsilon\in(0,1)$. Therefore,
Lemma \ref{lemE1} implies that $\P$-almost surely,
\begin{eqnarray}
\E^\pi_{\mathbf{s}} \biggl(\int_\Theta f(B(\mathbf{u}))\,
\d\mathbf{u} \biggr)&\ge& c^{-1}\mathrm{e}^{-c \|B(\mathbf{s})\|^2 - c
M^2}\nonumber\\[-8pt]\\[-8pt]
&\ge& c_{14}^{-1}\mathrm{e}^{-c_{14}\|B(\mathbf{s})\|^2}.\nonumber
\end{eqnarray}
On the other hand, Lemma \ref{lemE2} assures us that
%
\begin{equation}\qquad
\E^\pi_{\mathbf{s}} \biggl( \biggl|\int_\Theta f(B(\mathbf{u}))\,
\d\mathbf{u} \biggr|^2 \biggr)\le c_{15} \mathrm{e}^{c_{15}\|B(\mathbf{s})\|^2}
\cdot\mathrm{I}_{d-2N}(f) \qquad \mbox{a.s. }[\P].\vadjust{\goodbreak}
\end{equation}
We combine the preceding displays
together with the Paley--Zygmund inequality
and deduce that $\P$-almost surely,
\begin{eqnarray}\label{eqpfLB}
\P^\pi_{\mathbf{s}} \{ B(\mathbf{u})\in A^\epsilon\mbox{ for some
}\mathbf{u}
\in\Theta \} &\ge&
\P^\pi_{\mathbf{s}} \biggl\{ \int_\Theta f(B(\mathbf{u}))
\,\d\mathbf{u}>0 \biggr\}\nonumber\\[-8pt]\\[-8pt]
&\ge&\frac{\mathrm{e}^{-c_{16}\|B(\mathbf{s})\|^2}}{c_{16}\mathrm{I}_{d-2N}(f)}.\nonumber
\end{eqnarray}

Let $\mathcal{P}_{\mathrm{ac}}(A^{\epsilon/2})$ denote the
collection of all absolutely continuous
probability density functions that are supported on $A^{\epsilon/2}$.
It is the case that
%
\begin{equation}
\Bigl[\inf_{f\in\mathcal{P}_{\mathrm{ac}}(A^{\epsilon/2})}
\mathrm{I}_{d-2N}(f) \Bigr]^{-1}
\asymp{\rm Cap}_{d-2N}(A^{\epsilon/2}),
\end{equation}
where the implied constants depend only on $d$, $N$ and $M$
(\cite{Khbook}, Exercise~414, page~423). According to
Lemma \ref{lemRPT}, there exists an $\mathcal{F}_\pi(\mathbf{s})$-measurable
$\mu_\epsilon\in\mathcal{P}(A^{\epsilon/2})$ such that
%
\begin{equation}\label{eqeqmeas}
{\rm Cap}_{d-2N}(A^{\epsilon/2})
\asymp[\mathrm{I}_{d-2N}(\mu_\epsilon) ]^{-1},
\end{equation}
where the implied constants depend only on $d$, $N$ and $M$.
Let $\phi_\epsilon$ denote a smooth probability density function
with support in $B(0,\epsilon/2)=\{0\}^{\epsilon/2}$.
Then, $f=f_\epsilon:=\phi_\epsilon*\mu_\epsilon$ is
in $\mathcal{P}_{\mathrm{ac}}(A^\epsilon)$ and is
$\mathcal{F}_\pi(\mathbf{s})$-measurable.
We can apply~\eqref{eqpfLB} with this $f$, in order to obtain
the following: $\P$-almost surely,
%
\begin{equation}
\P^\pi_{\mathbf{s}} \{ B(\mathbf{u})\in A^\epsilon\mbox{ for some }
\mathbf{u}
\in\Theta\} \ge\frac{\mathrm{e}^{-c_{16}\|B(\mathbf{s})\|^2}}{c_{16}
\mathrm{I}_{d-2N}(\phi_\epsilon*\mu_\epsilon)}.
\end{equation}
But $\mathrm{I}_{d-2N}(\phi_\epsilon*\mu_\epsilon)\le C
\mathrm{I}_{d-2N}(\mu_\epsilon)$ for a finite nonrandom constant $C$ that
depends only on $d$, $N$, and $\sup\{|z|\dvtx z\in A\}$; see Theorems B.1
and B.2 of
\cite{DKN}. Therefore, we can deduce from \eqref{eqeqmeas} that
%
\begin{equation}\qquad
\P^\pi_{\mathbf{s}} \{ B(\mathbf{u})\in A^\epsilon\mbox{ for some }
\mathbf{u}\in\Theta \}
\ge c_{17} \mathrm{e}^{-c_{17}\|B(\mathbf{s})\|^2}
{\rm Cap}_{d-2N}(A^{\epsilon/2}).
\end{equation}
The resulting inequality
holds almost surely, simultaneously for all rational $\epsilon\in(0,1)$.
Therefore, we can let $\epsilon$ converge downward to zero, and
appeal to Choquet's capacitability theorem to deduce that
$\P$-almost surely,
%
\begin{equation}\qquad
\P^\pi_{\mathbf{s}} \{ B(\mathbf{u})\in A\mbox{ for some }\mathbf{u}
\in\Theta \}
\ge c_{18} \mathrm{e}^{-c_{17}\|B(\mathbf{s})\|^2}\cdot
\operatorname{Cap}_{d-2N}(A).
\end{equation}
(Choquet's theorem tells us that the preceding capacities are
outer regular; therefore as $A^\epsilon$ converges downward to
$A$, so do their respective capacities converge downward to the
capacity of $A$.)
Consequently,
%
\begin{equation}
\P^\pi_{\mathbf{s}} \{ B(\mathbf{u})\in A\mbox{ for some }\mathbf{u}
\in\Theta \}\unrhd\operatorname{Cap}_{d-2N}(A).
\end{equation}

We complete the theorem by deriving the converse direction;
that is,
%
\begin{equation}\label{eqGOAL1}
\P^\pi_{\mathbf{s}} \{ B(\mathbf{u})\in A\mbox{ for some }\mathbf{u}
\in\Theta \}\unlhd\operatorname{Cap}_{d-2N}(A).
\end{equation}
Equation \eqref{eqGOAL1} holds vacuously unless there is a positive probability
that the following happens:
%
\begin{equation}\label{eqassume}
\P^\pi_{\mathbf{s}} \{ B(\mathbf{u})\in A\mbox{ for some }\mathbf{u}
\in\Theta \}>0.
\end{equation}
Therefore, we may assume that \eqref{eqassume} holds with
positive probability without incurring
any further loss in generality.

Define
%
\begin{equation}\quad
T_1 := \inf\{ u_1\ge0\dvtx B(\mathbf{u})
\in A\mbox{ for some }\mathbf{u}=(u_1,\ldots,u_N)\in\Theta \},
\end{equation}
where $\inf\varnothing:=\infty$. Evidently $T_1$ is a
random variable with values in $\pi_1\Theta\cup\{\infty\}$,
where $\pi_l$ denotes the projection map which takes
$\mathbf{v}\in\R^N$ to $v_l$. Having constructed
$T_1,\ldots,T_j$ for $j\in\{1,\ldots,N-1\}$,
with values, respectively, in $\pi_1\Theta\cup\{\infty\}$,
\dots, $\pi_j\Theta\cup\{\infty\}$, we
define $T_{j+1}$ to be $+\infty$ almost surely on
$\bigcup_{l=1}^j\{T_l=\infty\}$, and
\[
T_{j+1} := \inf\{ u_{j+1}\ge0\dvtx
B(T_1,\ldots,T_j,u_{j+1},\ldots,u_N)\in A
\mbox{ for some }\mathbf{u}^T\in\Theta
\},
\]
almost surely on $\bigcap_{l=1}^j\{T_l<\infty\}$, where
in the preceding display
%
\begin{equation}
\mathbf{u}^T:= ( T_1,\ldots,T_j,u_{j+1},\ldots,u_N ).
\end{equation}
In this way, we obtain a random variable
$\mathbf{T}$, with values in $\Theta\cup\{\infty\}^N$,
defined as
%
\begin{equation}
\mathbf{T}:= ( T_1,\ldots,T_N ).
\end{equation}

Because \eqref{eqassume} holds with positive
probability, it follows that
%
\begin{equation}\label{eqkeyequiv}
\P^\pi_{\mathbf{s}} \{ \mathbf{T}\in\Theta \}
\asymp
\P^\pi_{\mathbf{s}} \{ B(\mathbf{u})\in A\mbox{ for some }\mathbf{u}
\in\Theta \}.
\end{equation}
If \eqref{eqassume} holds for some realization
$\omega\in\Omega$, then we define, for all
Borel sets $G\subseteq\R^d$,
%
\begin{equation}
\rho(G)(\omega) := \P^\pi_{\mathbf{s}} \bigl(
B(\mathbf{T})\in G  | \mathbf{T}\in\Theta \bigr)(\omega).
\end{equation}
Otherwise, we choose and fix some point $a\in A$
and define $\rho(G)(\omega):=\delta_a(G)$. It
follows that $\rho$ is a random $\F_\pi(\mathbf{s})$-measurable
probability measure on $A$.


Let $\phi_1 \in C^\infty(\R^d)$ be a probability density
function such that $\phi_1(x)=0$ if \mbox{$\|x\|>1$}. We define
an approximation to the identity $\{\phi_\eps\}_{\eps>0}$
by setting
%
\begin{equation} \label{Defphie}
\phi_\eps(x) := \frac{1}{\eps^d}\phi_1 \biggl( \frac{x}{\eps} \biggr)
\qquad \mbox{for all }x\in\R^d\mbox{ and }\eps>0.
\end{equation}
We plan to apply Lemma \ref{lemLB} with $f:=\rho*\psi_{\epsilon}$,
where $\psi_{\epsilon}(x) := \phi_{\epsilon/2}*\phi_{\epsilon/2}(x)$.
Furthermore, we can choose a good modification
of the conditional expectation in that lemma to deduce
that the null set off which the assertion fails can be
chosen independently of $\mathbf{s}$; see Lemma \ref{lemLBunif}.

Note that the support of $\rho*\psi_\epsilon$ is contained
in $A^\epsilon$. It follows from Lem\-ma~\ref{lemLB} that
$\P$-almost surely,
\begin{eqnarray}\label{eqkeyUB}
&&\sup_{\mathbf{t}\in\Theta}
\E^\pi_{\mathbf{t}} \biggl(\int_{\Theta^\eta}
(\rho*\psi_\epsilon)(B(\mathbf{u}))\,\d\mathbf{u}
\biggr)\nonumber\\[-9pt]\\[-9pt]
&&\qquad  \ge c_{11} \mathbf{1}_{\{\mathbf{T}\in\Theta\}}
\int_{\R^d}\kappa_{d-2N}(z) (\rho*\psi_\epsilon)
\bigl(z+B(\mathbf{T})\bigr)\,\d z .\nonumber
\end{eqnarray}
The constant $c_{11}$ is furnished by Lemma \ref{lemLB}.
Moreover, $\Theta^\eta$ denotes the closed $\eta$-enlargement
of $\Theta$.
We square both sides and take $\E^\pi_{\mathbf{s}}$-expectations.
Because $\mathbf{s}\prec_\pi\mathbf{t}$ for all $\mathbf{t}\in
\Theta$,
Lemma \ref{lemcairoli1} tells us that the $\E^\pi_{\mathbf{s}}$-expectation
of the square of the left-hand side of \eqref{eqkeyUB} is at most
%
\begin{equation}
4^N \sup_{\mathbf{t}\in\Theta}
\E^\pi_{\mathbf{s}} \biggl( \biggl| \E^\pi_{\mathbf{t}} \biggl[
\int_{\Theta^\eta} (\rho*\psi_\epsilon)(B(\mathbf{u}))\,
\d\mathbf{u} \biggr] \biggr|^2 \biggr).
\end{equation}
By the conditional form of Jensen's inequality,
$|\E^\pi_{\mathbf{t}} Z|^2\le\E^\pi_{\mathbf{t}}(Z^2)$ (a.s.),
for all square-integrable random variables $Z$. Moreover,
$\mathbf{s}\prec_\pi\mathbf{t}$ implies that
$\E^\pi_{\mathbf{s}}\E^\pi_{\mathbf{t}}=\E^\pi_{\mathbf{s}}$;
this follows from
the tower property of conditional expectations. Consequently,
%
\begin{eqnarray}
&&\E^\pi_{\mathbf{s}} \biggl( \biggl|\sup_{\mathbf{t}\in\Theta}
\E^\pi_{\mathbf{t}} \biggl(\int_{\Theta^\eta}
(\rho*\psi_\epsilon)(B(\mathbf{u}))\,\d\mathbf{u} \biggr) \biggr|^2 \biggr)\nonumber \\[-2pt]
&&\qquad  \le4^N\E^\pi_{\mathbf{s}} \biggl( \biggl|
\int_{\Theta^\eta} (\rho*\psi_\epsilon)(B(\mathbf{u}))\,
\d\mathbf{u} \biggr|^2 \biggr)\\[-2pt]
&&\qquad  \le c \mathrm{e}^{c\|B(\mathbf{s})\|^2} \cdot\mathrm{I}_{d-2N}(\rho*\psi_\epsilon),\nonumber
\end{eqnarray}
where the last inequality follows from Lemma \ref{lemE2}.
This and \eqref{eqkeyUB} together imply that with probability
one $[\P]$,
\begin{eqnarray}
&&c \mathrm{e}^{c\|B(\mathbf{s})\|^2} \cdot\mathrm{I}_{d-2N}(\rho
*\psi_\epsilon) \nonumber\\[-2pt]
&&\qquad \ge\E^\pi_{\mathbf{s}} \biggl( \biggl[ 
\int_{\R^d}\kappa_{d-2N}(z)
(\rho*\psi_\epsilon)
\bigl(z+B(\mathbf{T})\bigr)\,\d z\cdot\mathbf{1}_{\{\mathbf{T}\in\Theta\}}
\biggr]^2 \biggr) \nonumber\\[-9pt]\\[-9pt]
&&\qquad =\E^\pi_{\mathbf{s}} \biggl( \biggl[ 
\int_{\R^d}\kappa_{d-2N}(z)
(\rho*\psi_\epsilon)
\bigl(z+B(\mathbf{T})\bigr)\,\d z \biggr]^2
  \Big| \mathbf{T}\in\Theta\biggr) \nonumber\\[-2pt]
&&\qquad \quad {}\times\P^\pi_{\mathbf{s}}\{\mathbf{T}\in\Theta\}.\nonumber
\end{eqnarray}
%
We apply the Cauchy--Schwarz inequality and the definition
of $\rho$---in this order---to deduce from the preceding that
%
\begin{eqnarray}
&&c \mathrm{e}^{c\|B(\mathbf{s})\|^2} \cdot\mathrm{I}_{d-2N}
(\rho*\psi_\epsilon) \nonumber\hspace*{-35pt}\\[-2pt]
&&\qquad \ge \biggl[\E^\pi_{\mathbf{s}} \biggl(
\int_{\R^d}\kappa_{d-2N}(z) (\rho*\psi_\epsilon)
\bigl(z+B(\mathbf{T})\bigr)\,\d z
  \Big| \mathbf{T}\in\Theta\biggr) \biggr]^2
\,{\times}\,\P^\pi_{\mathbf{s}}\{\mathbf{T}\in\Theta\} \hspace*{-35pt}\\[-2pt]
&&\qquad = \biggl[\int_A \rho(\d x)\int_{\R^d} \,\d z\, \kappa_{d-2N}(z)
(\rho*\psi_\epsilon)(z+x)
\biggr]^2 \times\P^\pi_{\mathbf{s}}\{\mathbf{T}\in\Theta\}.\nonumber\hspace*{-35pt}
\end{eqnarray}
The term in square brackets is equal to
$\int(\kappa_{d-2N}*\rho*\psi_\epsilon)\,\d\rho$.
Since $\psi_\epsilon=\phi_{\epsilon/2}*\phi_{\epsilon/2}$,
that same term in square brackets is equal to
$\mathrm{I}_{d-2N}(\rho*\phi_{\epsilon/2})$. Thus, the following
holds $\P$-almost surely:
%
\begin{equation}\label{eqnearlydone}
c \mathrm{e}^{c\|B(\mathbf{s})\|^2} \cdot\mathrm{I}_{d-2N}(\rho*\phi
_{\epsilon/2}*
\phi_{\epsilon/2})
\ge [ \mathrm{I}_{d-2N}(\rho*\phi_{\epsilon/2})
]^2 \times\P^\pi_{\mathbf{s}}\{\mathbf{T}\in\Theta\}.\hspace*{-35pt}
\end{equation}

In order to finish the proof we now consider
separately the three cases where
$d<2N$, $d>2N$ and $d=2N$. If $d<2N$, then
\eqref{eqGOAL1} holds because the right-hand side is 1.

If $d>2N$, then Theorem B.1 of Dalang et al. \cite{DKN}
tells us that
%
\begin{equation}
\mathrm{I}_{d-2N}(\rho*\phi_{\epsilon/2}*\phi_{\epsilon/2})
\le\mathrm{I}_{d-2N}(\rho*\phi_{\epsilon/2}).
\end{equation}
Since $\kappa_{d-2N}$ is lower semicontinuous, Fatou's
lemma shows that
%
\begin{equation}
\liminf_{\epsilon\downarrow0} \mathrm{I}_{d-2N}(\rho*\phi_{
\epsilon/2})\ge\mathrm{I}_{d-2N}(\rho).
\end{equation}
Therefore, \eqref{eqnearlydone} implies that:
\begin{itemize}[(ii)]
\item[(i)] $\mathrm{I}_{d-2N}(\rho)<\infty$ [thanks also to \eqref
{eqassume}];
\item[(ii)] $\P^\pi_{\mathbf{s}}\{ \mathbf{T}\in\Theta\}\le c
\exp(c\|B(\mathbf{s})\|^2)/
\mathrm{I}_{d-2N}(\rho)$ almost surely.
\end{itemize}
This proves that $\P$-almost surely,
\begin{eqnarray}\label{eqGG}
\P^\pi_{\mathbf{s}} \{ \mathbf{T}\in\Theta\}
&\le&\frac{c
\mathrm{e}^{c\|B(\mathbf{s})\|^2}}{\mathrm{I}_{d-2N}(\rho)}\nonumber\\[-8pt]\\[-8pt]
&\le& c \mathrm{e}^{c\|B(\mathbf{s})\|^2} \cdot\operatorname{Cap}_{d-2N}(A).\nonumber
\end{eqnarray}
Consequently, \eqref{eqkeyequiv} implies the theorem
in the case that $d>2N$.

The final case that $d=2N$ is
handled similarly, but this time we use Theorem~B.2 of
Dalang et al. \cite{DKN} in place of their Theorem B.1.
\qed

\section{\texorpdfstring{Proofs of Theorem \protect\ref{thmain} and its corollaries}{Proofs of Theorem 1.1 and its corollaries}}\label{secpfmain}

We start with the following result which deals with intersections of
the images of the Brownian sheet of disjoint boxes that satisfy
certain configuration conditions.

\begin{theorem}\label{Th51}
Let $\Theta_1,\ldots,\Theta_k$ in $(0,\infty)^N$ be disjoint,
closed and
nonrandom upright boxes that satisfy the following properties:
\begin{enumerate}[(0)]
\item[(1)]
for all $j=1,\ldots,k-1$ there exists $\pi(j)\subseteq\{1,\ldots,N\}$
such that $\mathbf{u}\prec_{\pi(j)}\mathbf{v}$
for all $\mathbf{u}\in\bigcup_{l=1}^j\Theta_l$ and $\mathbf{v}\in
\Theta_{j+1}$;
\item[(2)] there exists a nonrandom $\eta>0$ such that
$|\mathbf{r}-\mathbf{q}|_\infty\ge\eta$ for all $\mathbf{r}\in
\Theta_i$
and $\mathbf{q}\in\Theta_j$, where $1\le i\neq j\le k$.
\end{enumerate}
Then for any Borel set $A \subseteq\R^d$,
%
\begin{equation}\label{EqEquiv}
 \P\Biggl\{ \bigcap_{j=1}^k B(\Theta_j) \cap A \ne\varnothing \Biggr\}>0
\quad  \Longleftrightarrow\quad
\P\Biggl\{ \bigcap_{j=1}^k W_j(\Theta_j) \cap A \ne\varnothing\Biggr\}>0,\hspace*{-35pt}
\end{equation}
where $W_1,\ldots,W_k$ are $k$ independent
$N$-parameter Brownian sheets in $\R^d$ (which are unrelated to $B$).
\end{theorem}

\begin{pf}
Under assumptions (1) and (2),
we can choose and fix nonrandom time points $\mathbf{s}_1,
\ldots,\mathbf{s}_{k-1}\in(0,\infty)^N$
such that for all $l=1,\ldots,k-1$:
\begin{enumerate}[(3)]
\item[(3)] $\mathbf{s}_l\prec_{\pi(l)}\mathbf{v}$
for all $\mathbf{v}\in\Theta_{l+1}$;
\item[(4)] $\mathbf{s}_l\succ_{\pi(l)}\mathbf{u}$
for all $\mathbf{u}\in\bigcup_{j=1}^l\Theta_j$.
\end{enumerate}
Because the elements of $\Theta_k$ dominate those of
$\Theta_1,\ldots,\Theta_{k-1}$ in partial order \mbox{$\pi(k-1)$},
Theorem \ref{thCPT} can be applied (under
$\P^{\pi(k-1)}_{\mathbf{s}_{k-1}}$) to show that for all nonrandom
Borel sets $A\subset\R^d$,
%
\begin{equation}\quad
\P\{ [\mathbf{B}]_k \cap A \neq\varnothing\}>0
 \quad \Longleftrightarrow\quad  \E\bigl[\operatorname{Cap}_{d-2N}
([\mathbf{B}]_{k-1} \cap A ) \bigr]>0,
\end{equation}
where
%
\begin{equation}
[\mathbf{B}]_k:= \bigcap_{j=1}^k B(\Theta_j).
\end{equation}
The main result of Khoshnevisan and Shi \cite{KS} is that
$\operatorname{Cap}_{d-2N}(E)>0$ is necessary and sufficient
for $\P\{W_k(\Theta_k)\cap E\neq\varnothing\}$ to be (strictly) positive,
where $W_k$ is a Brownian sheet that is independent of
$B$. We apply this with $E:= [\mathbf{B}]_{k-1} \cap A$ to deduce that
%
\begin{equation}\qquad
\P\{ [\mathbf{B}]_k \cap A\neq\varnothing\}>0
 \quad \Longleftrightarrow\quad  \P\{
[\mathbf{B}]_{k-1} \cap W_k(\Theta_k) \cap A\neq
\varnothing\} >0.
\end{equation}
Because $W_k$ is independent of
$B$, and thanks to (3) and (4) above,
we may apply Theorem \ref{thCPT} inductively to deduce that
%
\begin{equation}\label{eqIE}\quad
\P\{ [\mathbf{B}]_k \cap A\neq\varnothing\}>0
\quad \Longleftrightarrow\quad  \P\Biggl\{
\bigcap_{j=1}^k W_j(\Theta_j)\cap A\neq
\varnothing\Biggr\} >0,
\end{equation}
where $W_1,\ldots,W_k$ are i.i.d. Brownian sheets. This proves
Theorem \ref{Th51}.
\end{pf}

Note that conditions (1) and (2) in Theorem \ref{Th51} are
satisfied for $k=2$ for two arbitrary upright boxes $\Theta_1$ and
$\Theta_2$ that have disjoint projections on each coordinate
hyperplane $s_i = 0$, $i=1,\dots,N$. Hence we are ready to derive
Theorem 1.1.

\begin{pf*}{Proof of Theorem \ref{thmain}}
Observe that there exist distinct points $\mathbf{s}$ and $\mathbf{t}
\in
(0,\infty)^N$ with $B(\mathbf{s})=B(\mathbf{t})\in A$, and such that
$s_i\neq t_i$, for all $i=1,\dots,N$, if and only if we can
find disjoint closed upright boxes $\Theta_1$ and $\Theta_2$, with
vertices with rational coordinates, such that $[\mathbf{B}]_2\cap
A\neq\varnothing$. Moreover, we may require $\Theta_1$ and $\Theta
_2$ to be such that
assumptions (1) and (2) of Theorem~\ref{Th51}
are satisfied. Since the family of pairs of such closed upright
boxes~$\Theta_1$ and~$\Theta_2$ is countable, it follows
that~\eqref{EqEquiv} implies Theorem~\ref{thmain}.
\end{pf*}

In order to apply Theorem \ref{thmain} to study the nonexistence
of double points of the Brownian sheet, we first provide some
preliminary results on the following subset of $M_2$:
\[
M_2^{(1)}:= \left\{x\in\R^d \left|
\matrix{B(\mathbf{s}_1)=B(\mathbf{s}_2) =x \mbox{ for
distinct }\mathbf{s}_1, \mathbf{s}_2\in(0,\infty)^N \cr
\mbox{with at least one common coordinate}\hspace*{46pt}
}\right.\right\}.
\]

Note that $M_2^{(1)}$ cannot be studied by using Theorem
\ref{thmain}. The next lemma will help us to show that
$M_2^{(1)}$ has negligible effect on the properties of $M_2$.

\begin{lemma}\label{LemExtra}
The random set $M_2^{(1)}$ has the following properties:
\begin{enumerate}[(ii)]
\item[(i)] $\dimh M_2^{(1)} \le4N - 2 -d$ a.s., and ``$\dimh M_2^{(1)}<0$''
means ``$M_2^{(1)}=\varnothing$'';

\item[(ii)] for every nonrandom Borel set $A\subseteq\R^d$,
%
\begin{equation}\label{EqM2hit}
\P\bigl\{M_2^{(1)} \cap A \ne\varnothing\bigr\} \le
\mathrm{const}\cdot\mathcal{H}_{2d- 2(2N-1)}(A),
\end{equation}
where $\mathcal{H}_{\beta}$ denotes the $\beta$-dimensional Hausdorff
measure.
\end{enumerate}
\end{lemma}

\begin{pf}
Part (i) follows from (ii) and a standard covering argument;
see, for example, \cite{BLX,DKN,Xiao09}. We omit the details and only
give the following rough outline. We only consider the case where
$\mathbf{s}_1,
\mathbf{s}_2 \in(0,\infty)^N$ are distinct, but have the same first
coordinates. This causes little loss of generality.

For a point in a fixed unit cube of $\R^N$, say $[1,2]^N$, there are $2^{2n}$
possible first coordinates of the form $1+i2^{-2n}$, $i=0,\dots,2^{2n} -1$.

For any given such first coordinate, there are $(2^{2n})^{N-1}$ points
in $[1,2]^N$
with all other coordinates of the same form as the first coordinate. In another
unit cube, such as $[1,2] \times[3,4]^{N-1}$, there are also $(2^{2n})^{N-1}$
points with a given first coordinate and all other coordinates of the form
$1+i2^{-2n}$, $i=0,\dots,2^{2n} -1$.

We cover the set $M_2^{(1)} \cap[0,1]^d$ by small boxes with sides of
length $n 2^{-n}$.
If we cover $[0,1]^d$ by a grid of small boxes with sides of length $n
2^{-n}$, the
probability that any small box $C$ in $[0,1]^d$ is needed to help cover
$M_2^{(1)}
\cap[0,1]^d$ because of the behavior of $B$ near $(u_1,u_2)$ and
$(u_1,v_2)$ is approximately
%
\begin{equation}\quad
\P\{B(u_1,u_2) \in C, \Vert B(u_1,u_2) - B(u_1,v_2)\Vert\leq n
2^{-n}\} \simeq(n 2^{-n})^{2d},
\end{equation}
where $(u_1,u_2) \in[1,2]^N$ and $(u_1,v_2) \in[1,2] \times
[3,4]^{N-1}$. Therefore, for $\gamma>0$,
%
\begin{equation}\label{rda1}\quad
 \E\Bigl( \sum(n 2^{-n})^\gamma\Bigr) \simeq2^{nd} (n 2^{-n})^\gamma\P\{
\mbox{a given small box is in the covering} \},\hspace*{-35pt}
\end{equation}
where the sum on the left-hand side is over all small boxes in a
covering of $M_2^{(1)} \cap[0,1]^d$.
The probability on the right-hand side is approximately
\begin{eqnarray}
&& \#\{\mbox{points } (u_1,u_2)\mbox{ and } (v_1,v_2) \mbox{ to be
considered} \} (n 2^{-n})^{2d}\nonumber\\[-8pt]\\[-8pt]
&&\qquad  = 2^{2n} \bigl((2^{2n})^{(N-1)} \bigr)^2 (n 2^{-n})^{2d}.\nonumber
\end{eqnarray}
It follows that the left-hand side of \eqref{rda1} is approximately
equal to
%
\begin{equation}
n^{\gamma+2d} (2^{-n})^{\gamma-4N+2+d}.
\end{equation}
This converges to $0$ if $\gamma> 4N - 2 - d$, and this explains
statement (i).

In order to prove (ii), we start with a hitting probability
estimate for $M_2$. Let $D := D_1\times D_2\times D_3$ denote a
compact upright box in $ (0, \infty)^{1+2(N-1)}$, where $D_2,
D_3 \subset(0, \infty)^{N-1}$ are disjoint. By using the argument
in the proof of Proposition 2.1 in Xiao \cite{Xiao99} we can show
that simultaneously for all $(a_1, \mathbf{a}_2, \mathbf{a}_3) \in
D$, $r > 0$ and $x
\in\R^d$,
%
\begin{equation}\label{Eqdouble-hit}
 \P\left\{
\matrix{
\exists (t_1,\mathbf{t}_2,\mathbf{t}_3)\in(a_1-r^2, a_1+r^2)
\times U_r(\mathbf{a}_2)\times U_r(\mathbf{a}_3)\cr
\mbox{such that } |B(t_1, \mathbf{t}_2)-x | \le r, |B(t_1,
\mathbf{t}_3)-x | \le r\hspace*{20pt}
}
\right\} =O(r^{2d}),\hspace*{-35pt}
\end{equation}
as $r\downarrow0$, where $U_r(\mathbf{a}):= \{\mathbf{t} \in\R
^{N-1}\dvtx |\mathbf{t}-
\mathbf{a}| \le r^2\}$. The proof of \eqref{Eqdouble-hit} is somewhat
lengthy. Since it is more or less a standard proof, we omit the details,
and offer instead only the following rough outline:
(a) For fixed $(t_1,\mathbf{t}_2,\mathbf{t}_3)$ in
$(a_1-r^2,a_2+r^2)\times
U_r(\mathbf{a}_2)\times U_r(\mathbf{a}_3)$, we have
$\P\{|B(t_1,\mathbf{t}_2)|\le r\}=O(r^d)$ thanks to direct computation;
(b) $\P\{|B(t_1,\mathbf{t}_3)|\le r \mid|B(t_1,\mathbf{t}_2)|$ $\le
r\} =O(r^d)$
because $B(t_1,\mathbf{t}_2)$ and $B(t_1,\mathbf{t}_3)$ are ``sufficiently
independent''; (c) $B$ in time-intervals of side length $r^2$ is roughly
``constant'' to within at most $r$ units.
Part (ii) follows from \eqref{Eqdouble-hit} and another covering argument
\cite{BLX,DKN}.
\end{pf}

We now show how Theorem \ref{thmain} can be
combined with the elegant theory of Peres \cite{Peres}
and Lemma \ref{LemExtra} to imply the corollaries
mentioned in the \hyperref[sec1]{Introduction}.


\begin{pf*}{Proof of Corollary \ref{cohit}}
Theorem 1.1 of Khoshnevisan and Shi \cite{KS} asserts that for each
$\nu\in\{1,2\}$, all nonrandom Borel sets $A\subset\R^d$, contained
in a
fixed compact subset of $\R^d$, and all upright boxes
$\Theta:=\prod_{j=1}^N[a_j,b_j]\subset(0,\infty)^N$,
there is a finite constant $R\ge1$ such that
%
\begin{equation}\label{eqKSPS}\quad
R^{-1} \operatorname{Cap}_{d-2N}(A) \leq
\P\{ W_\nu(\Theta)\cap A\neq\varnothing\}
\leq R\operatorname{Cap}_{d-2N}(A).
\end{equation}

We first consider the case where $d>2N$. Because $W_1$ and $W_2$ are
independent,
Corollary 15.4 of Peres (\cite{Peres}, page 240) and (\ref{eqKSPS})
imply that for
all upright boxes $\Theta_1,\Theta_2 \subset(0, \infty)^N$
%
\begin{equation}\label{eqLambda1}
 \P\{ W_1(\Theta_1)\cap
W_2(\Theta_2)\cap A\neq\varnothing\}>0
\quad \Longleftrightarrow\quad
\operatorname{Cap}_{2(d-2N)}(A)>0.\hspace*{-35pt}
\end{equation}

Next, let us assume that $\operatorname{Cap}_{2(d-2N)}(A)>0$.
We choose arbitrary upright boxes $\Theta_1$
and $\Theta_2$ that have disjoint projections on each coordinate
hyperplane $s_i = 0$,
$i=1,\dots,N$. It follows that
$\P\{M_2 \cap A \ne\varnothing \} > 0$, thanks to~\eqref{eqLambda1} and Theorem~\ref{thmain}.

In order to prove the converse, we assume that $\operatorname
{Cap}_{2(d-2N)}(A)=0$. Then $\dim A
\le2(d-2N)$ which implies $\mathcal{H}_{2d- 2(2N-1)}(A) =0$.
It follows from Lem\-ma~\ref{LemExtra} that $\P\{M_2^{(1)} \cap A \ne
\varnothing \} = 0.$ On the other hand, \eqref{eqLambda1} and Theorem~\ref{thmain} imply
that $\P\{(M_2\setminus M_2^{(1)}) \cap A \ne
\varnothing \} = 0.$ This finishes the proof when $d > 2N$.

In the case $d=2N$, where~\eqref{eqbilog} appears,
we also use~\eqref{eqKSPS}, Corollary~15.4 of Peres \cite{Peres}, page
240,
and Lemma~\ref{LemExtra}.

Finally, if $2N > d$, then $B$ hits points by \eqref{eqKSPS}.
This implies the last conclusion in Corollary \ref{cohit}.
\end{pf*}

\begin{pf*}{Proof of Corollary \ref{coMP}}
We appeal to Corollary \ref{cohit} with $A:=\R^d$,
and use \eqref{eqtaylor} to deduce that
$\P\{M_2\neq\varnothing\}>0$ if and only if $2(d-2N)<d$. Next, we
derive the second assertion of the corollary (Fristedt's conjecture).

Choose and fix some $x\in\R^d$. Corollary \ref{cohit}
tells us that
\mbox{$\P\{x\in M_2\}>0$} if and only if $\operatorname{Cap}_{2(d-2N)}(\{x\})>0$.
Because the only probability measure on~$\{x\}$
is the point mass, the latter positive-capacity
condition is equivalent to the condition that $d<2N$. According to the Tonelli
theorem, $\E(\operatorname{meas} M_2)=\int_{\R^d}\P\{x\in M_2\}\,
\d x$,
where ``$\operatorname{meas}M_2$'' denotes
the $d$-dimensio\-nal
Lebesgue measure of $M_2$. It follows
readily from this discussion that
%
\begin{equation}\label{eqEmeas}
\E(\operatorname{meas} M_2)>0 \quad \Longleftrightarrow\quad
d<2N.
\end{equation}
If $d\ge2N$, then this proves that $\operatorname{meas} M_2=0$ almost surely.

It only remains to show that if $d<2N$, then $\operatorname{meas} M_2
>0$ almost surely.
For any integer $\ell\ge0$, define
\[
M_{2,\ell} := \{ x\in\R^d\dvtx B(\mathbf{s}^1)=
B(\mathbf{s}^2)=x\mbox{ for distinct }\mathbf{s}^1,
\mathbf{s}^2\in[ 2^{\ell},2^{\ell+1} ]^N \}.
\]
Given a fixed point $x\in\R^d$, the scaling properties of the Brownian
sheet imply that
%
\begin{equation}\label{eqscaling}
\P\{x\in M_{2,\ell}\}=\P\{2^{-\ell N/2}x\in M_{2,0}\}.
\end{equation}
By using Theorem \ref{thmain} and \eqref{eqLambda1}, we see that
%
\begin{equation}
\inf_{x\in[-q,q]^d}\P\{x\in M_{2,0}\}>0
\qquad \mbox{for all }q>0.
\end{equation}
The scaling property \eqref{eqscaling} then implies that
%
\begin{equation}
\gamma_q:=\inf_{\ell\ge0}
\inf_{x\in[-q,q]^d}\P\{x\in M_{2,\ell}\}>0\qquad
\mbox{for all }q>0.
\end{equation}
In particular,
%
\begin{equation}\label{OP}
\P\{ x\in M_{2,\ell}
\mbox{ for infinitely many }\ell\ge0 \}
\ge\gamma_q
>0 \quad \mbox{for all }x\in[-q,q]^d.\hspace*{-40pt}
\end{equation}
By the zero--one law of Orey and Pruitt
(\cite{OP}, pages 140--141), the left-hand
side of \eqref{OP} is identically
equal to one. But that left-hand side is at most
$\P\{x\in M_2\}$. Because
$q$ is arbitrary, this proves that $\P\{x\in M_2\}=1$ for all $x\in\R^d$
when $d<2N$. By Tonelli's theorem,
$\P\{ x\in M_2 \mbox{ for almost all $x\in\R^d$} \}=1$, whence
$\operatorname{meas} M_2=\infty$ almost surely,
and in particular $\operatorname{meas} M_2>0$ almost surely.
\end{pf*}

\begin{pf*}{Proof of Proposition \ref{comonotone}}
According to \eqref{eqKSPS} and Corollary 15.4 of Peres \cite{Peres}, page
240,
the following is valid for all $k$ upright boxes $\Theta_1,\ldots
,\Theta_k\subset
(0,\infty)^N$ with vertices with rational coordinates:
%
\begin{equation}\label{rd516}
\P\Biggl(\bigcap_{\nu= 1}^k W_\nu(\Theta_\nu) \cap A \neq
\varnothing \Biggr) >0 \quad \Longleftrightarrow\quad  \operatorname{Cap}_{k(d-2N)}(A)>0.
\end{equation}

Observe that $\P\{ \widetilde{M}_k \cap A \neq\varnothing\} >0$ if and
only if there exists a partial order $\pi\subseteq\{1,\dots,N\}$
together with $k$ disjoint upright boxes $\Theta_1,\dots,\Theta_k$ in
$(0,\infty)^N$, with vertices with rational coordinates, such that for
$1\leq i < j \leq k$, $\mathbf{s} \in\Theta_i$ and $\mathbf{t} \in
\Theta_j$
implies $\mathbf{s} \ll_\pi\mathbf{t}$, and
%
\begin{equation}
\P\Biggl(\bigcap_{\nu= 1}^k B(\Theta_\nu) \cap A \neq\varnothing \Biggr) >0.
\end{equation}
In addition, $\Theta_1,\dots,\Theta_k$ can be chosen so as to
satisfy (1) and (2) of Theorem~\ref{Th51} [with $\pi(j) = \pi$,
$j=1,\dots,k-1$].
It follows from Theorem~\ref{Th51} that
%
\begin{equation}\qquad
\P\{ \widetilde M_k \cap A \neq\varnothing\} >0
\quad \Longleftrightarrow\quad
\P\Biggl(\bigcap_{\nu= 1}^k W_\nu(\Theta_\nu) \cap A \neq\varnothing
\Biggr) >0.
\end{equation}
Owing to \eqref{rd516}, the right-hand side is equivalent to
the [strict] positivity of $\operatorname{Cap}_{k(d-2N)}(A)$;\vspace*{2pt}
this proves the first statement in Proposition \ref{comonotone}, and
the second statement follows by taking \eqref{eqtaylor} into account.
\end{pf*}

\begin{remark}
The following is a consequence of Proposition \ref{comonotone}:
Fix an integer $k > 2$, and suppose that with positive probability
there exist distinct $\mathbf{u}_1,\ldots,\mathbf{u}_k\in
(0,\infty)^N$ such that $W_1(\mathbf{u}_1)=\cdots=
W_k(\mathbf{u}_k)\in A$. Then with positive probability
there exist distinct $\mathbf{u}_1,\ldots,
\mathbf{u}_k\in(0,\infty)^N$ such that $B(\mathbf{u}_1)=
\cdots=B(\mathbf{u}_k)\in A$. We believe the converse is true.
But Proposition~\ref{comonotone}, and
even~\eqref{eqIE}, implies the converse only for special configurations
of $\mathbf{u}_1, \ldots,\mathbf{u}_k$.
In particular, the question of the existence of
$k$-multiple points in critical dimensions [$k>2$ for which
$k(d-2N)=d$] remains open.
\qed
\end{remark}

\begin{pf*}{Proof of Corollary \ref{comain}}
We can combine \eqref{EqEquiv} and \eqref{eqLambda1} with
Corollary~\ref{cohit} and deduce\vadjust{\goodbreak} that whenever $\Theta_1$
and $\Theta_2$ are the upright boxes of the proof of
Theorem \ref{thmain},
%
\begin{equation}\label{eqhitB}\qquad
\P\{ B(\Theta_1)\cap B(\Theta_2)
\cap A \neq\varnothing\}>0
\quad \Longleftrightarrow\quad
\operatorname{Cap}_{2(d-2N)}(A)>0.
\end{equation}
This is valid for all nonrandom Borel sets $A\subseteq\R^d$.

By Frostman's theorem, if $\dimh A<2(d-2N)$, then
$\operatorname{Cap}_{2(d-2N)}(A)=0$; see \eqref{eqfrostman}. Consequently,
Corollary \ref{cohit}
implies that $M_2\cap A=\varnothing$ almost surely.

Next, consider the case where $\dimh A\ge2(d-2N)>0$. Choose
and fix some constant $\rho\in(0,d)$.
According to Theorem 15.2 and Corollary 15.3 of
Peres \cite{Peres}, pages 239--240,
we can find a random set $\mathbf{X}_\rho$, independent of the
Brownian sheet~$B$, that has the following properties:
\begin{itemize}[--]
\item[--] for all nonrandom Borel sets
$A\subseteq\R^d$,
%
\begin{equation}\label{eqP1}
\P\{\mathbf{X}_\rho\cap A\neq\varnothing\}>0
\quad \Longleftrightarrow\quad
\operatorname{Cap}_\rho(A) >0;
\end{equation}
\item[--] for all nonrandom Borel sets
$A\subseteq\R^d$ and all $\beta>0$,
%
\begin{equation}\label{eqP2}
\P\{ \operatorname{Cap}_\beta(\mathbf{X}_\rho\cap A)>0 \}>0
\quad \Longleftrightarrow\quad
\operatorname{Cap}_{\rho+\beta}(A)>0.
\end{equation}
\end{itemize}
[Indeed, $\mathbf{X}_\rho$ is the fractal-percolation
set $Q_d(\kappa_\rho)$ of Peres (loc. cit.).]

Equation \eqref{eqhitB} implies that
%
\begin{equation}\label{eqG1}\qquad
\P\{ [\mathbf{B}]_2
\cap\mathbf{X}_\rho\cap A \neq\varnothing  | \mathbf{X}_\rho\}>0
\quad \Longleftrightarrow\quad
\operatorname{Cap}_{2(d-2N)}(\mathbf{X}_\rho\cap A)>0,
\end{equation}
where we recall that $[\mathbf{B}]_2:=B(\Theta_1)\cap B(\Theta_2)$.
Thanks to \eqref{eqP2},
%
\begin{equation}\label{eqG2}
\P\{ [\mathbf{B}]_2
\cap\mathbf{X}_\rho\cap A \neq\varnothing  | \mathbf{X}_\rho\}
\asymp\operatorname{Cap}_{2(d-2N)+\rho}(A)
\end{equation}
holds almost surely. At the same time, \eqref{eqP1} implies that
%
\begin{equation}
\P\{ [\mathbf{B}]_2
\cap\mathbf{X}_\rho\cap A \neq\varnothing
| B \} \asymp
\operatorname{Cap}_{\rho} ( [\mathbf{B}]_2\cap
A ).
\end{equation}
Therefore, we compare the last two displays to deduce that
%
\begin{equation}\label{eqG3}
\P\{ \operatorname{Cap}_{\rho} ( [\mathbf{B}]_2 \cap
A ) >0 \} > 0 \quad \Longleftrightarrow\quad
\operatorname{Cap}_{2(d-2N)+\rho}(A) >0.
\end{equation}
Frostman's theorem (\cite{Khbook}, page 521) then implies
the following:
%
\begin{equation}\label{eqG4}
\| \dimh( [\mathbf{B}]_2 \cap
A ) \|_{L^\infty(\P)} = \dimh A - 2(d-2N).
\end{equation}
This and (\ref{EqM2Adim}) together imply readily the announced formula
for the $\P$-essential supremum of $\dimh(M_2\cap A)$.

The remaining case is when $d=2N$. In that case,
we define for all measurable
functions $\kappa\dvtx\R_+\to\R_+\cup\{\infty\}$,
%
\begin{equation}
\operatorname{Cap}_\kappa(A) := [ \inf\mathrm{I}_\kappa(\mu) ]^{-1},
\end{equation}
where the infimum is taken over all compactly supported
probability measures $\mu$ on $A$,
and
%
\begin{equation}
\mathrm{I}_\kappa(\mu) := \iint\kappa(\|x-y\|) \mu(\d x)
\mu(dy).
\end{equation}
Then the preceding argument
goes through, except we replace:
\begin{itemize}[--]
\item[--] $\operatorname{Cap}_{2(d-2N)}(\mathbf{X}_\rho\cap A)$
by $\operatorname{Cap}_f(\mathbf{X}_\rho\cap A)$ in \eqref{eqG1},
where $f(u):=|\log_+(1/u)|^2$;
\item[--] $\operatorname{Cap}_{2(d-2N)+\rho}(A)$
by $\operatorname{Cap}_g(A)$ in \eqref{eqG2} and \eqref{eqG3}, where
$g(u) := |u|^\rho f(u)$;
\item[--] $2(d-2N)$ by zero
on the right-hand side of \eqref{eqG4}.
\end{itemize}
The justification for these replacements is the same as
for their analogous assertions in the case $d>2N$.
This completes our proof.
\end{pf*}


%

\printaddresses

\end{document}